\documentclass[a4paper]{article}
\usepackage{latexsym}
\usepackage{amsthm,amsmath,amssymb}
\date{}

\usepackage{accents}

\newtheorem{thm}{Theorem}

\newtheorem{prop}[thm]{Proposition}

\newtheorem{ex}[]{Example}

\newtheorem{cor}[thm]{Corollary}

\newtheorem{lem}[thm]{Lemma}
\newtheorem{sublem}[thm]{Sublemma} 
\newtheorem{defn}[thm]{Definition}

\newcommand{\Q}{\mathbb{Q}}
\def\O{\mathcal O}

\def\Z{\mathbb Z} \def\R{\mathbb R}  
\def\={\;=\;}  \def\+{\,+\,}  \def\Q{\Bbb Q}  \def\Z{\Bbb Z}  \def\F{\Bbb F} \def\FF{\mathcal F}    
\def\L{\mathcal L} 
\def\A{\mathbb A}

    \def\m{\bold m}
    
\def\cts{\mathrm{cts}}

\def\be{\begin{equation}}   \def\ee{\end{equation}}
\def\bes{\begin{equation*}}   \def\ees{\end{equation*}}

\newcommand{\ar}{\mathrm{ar}}
\newcommand{\coef}{\mathrm{coeft}}

\newcommand{\Spec}{\mathrm{Spec}}
\newcommand{\fin}{\mathrm{fin}}
\newcommand{\Frac}{\mathrm{Frac}}

\newcommand{\Res}{\mathrm{Res}}
\newcommand{\res}{\mathrm{res}}

\newcommand{\ord}{\mathrm{ord}}

\newcommand{\red}{\mathrm{red}}
\newcommand{\sep}{\mathrm{sep}}
\newcommand{\Tr}{\mathrm{Tr}}

\newcommand{\law}{\longleftarrow}
\newcommand{\raw}{\longrightarrow}

\newcommand{\dis}{\displaystyle}

\begin{document}
\title{\bf Arithmetic Cohomology Groups}
\author{{\bf K. Sugahara}\ \  {\small{and}}\ \ \  {\bf L. Weng}} 
\maketitle
\begin{abstract} We first introduce global arithmetic cohomology groups for quasi-coherent sheaves on arithmetic 
varieties, adopting an adelic approach. Then, we establish fundamental properties, such as topological duality and
inductive long exact sequences, for these cohomology groups on arithmetic surfaces. Finally, we expose basic 
structures for ind-pro topologies on adelic spaces of arithmetic surfaces. In particular, we show that these adelic 
spaces are topologically self-dual.
\end{abstract}
\begin{tableofcontents}
\end{tableofcontents}
\section*{Introduction}

In the study of arithmetic varieties, cohomology theory has been developed along with the line of establishing an
intrinsic relation between arithmetic Euler characteristics and arithmetic intersections. For examples, for an arithmetic 
curve $\Spec\,\O_F$ associated to the integer ring $\O_F$ of a number field $F$ with discriminant $\Delta_F$ and 
a metrized vector sheaf $\overline E$ on it, we have the Arakelov-Riemann-Roch formula
$$
\chi_{\ar}(F, \overline E)=\deg_\ar(\overline E)-\frac{\mathrm{rank}E}{2}\log|\Delta_F|;
$$
And, for a regular arithmetic surface $\pi:X\to \Spec\,\O_F$ and a metrized line sheaf $\overline \L$ on it, if we equip 
with $X_\infty$ a K\"ahler metric, and  line sheaves $\lambda(\L)$ and $\lambda(\O_X)$ with the Quillen metrics,  
namely, equip determinants of relative cohomology groups  with determinants of $L^2$-metrics modified by analytic 
torsions, then we have the Faltings-Deligne-Riemann-Roch  isometry
$$
\overline {\lambda(\L)}^{\otimes 2}\otimes\overline{\lambda(\O_X)}^{\otimes -2}
\simeq
\langle\overline\L,(\overline \L\otimes \overline K_\pi^{\otimes -1})\rangle;
$$
More generally, for higher dimensional arithmetic varieties, we have the works of (Bismut-)Gillet-Soul\'e.

In this paper, we start to develop a genuine cohomology theory for arithmetic varieties, as a continuation of  the works 
of Parshin ([P1,2]), Beilinson ([B]), Osipov-Parshin ([OP]), and  our own study ([W]). Our aims here are to construct 
arithmetic cohomology groups $H_{\ar}^i$  for quasi-coherent sheaves on arithmetic varieties, and to establish 
topological dualities among these cohomology groups for arithmetic surfaces. 

The approach we take in this paper is an adelic one. Here, we use three main ideas form the literatures. Namely, the 
first one of adelic complexes initiated in the classical works [P1,2] and [B], (see also [H]), which is recalled in \S1.1 and 
used in \S1.2 systematically; the second  one of ind-pro structures over adelic spaces from [OP], which is recalled 
in \S1.2.1 and motivates our general constructions in \S1.2.2; and the final one on the uniformity structure between 
finite and infinite adeles from [W], which is recalled in \S1.2.4 and plays an essential role in \S1.2.3 when we construct 
our adelic spaces. In particular, we are able to introduce arithmetic adelic complexes $(\A_\ar^*(X,\FF),d^*)$ for 
quasi-coherent sheaves $\FF$ over arithmetic varieties $X$, and hence are able to define their associated arithmetic 
cohomology groups $H_\ar^i(X,\FF):=H^i(\A_\ar^*(X,\FF),d^*)$. Consequently, we have the following

\vskip 0.20cm
\noindent
{\bf Theorem I.} 
{\it Let $X$ be an arithmetic variety and $\FF$ be a quasi-coherent sheaf on $X$, then there exist a natural arithmetic 
adelic complex $(\A_\ar^*(X,\FF),d^*)$ and hence arithmetic cohomology groups 
$H_\ar^i(X,\FF):=H^i(\A_\ar^*(X,\FF),d^*)$. In particular, $H^i_{\ar}(X,\FF)=0$ unless $i=0,1,\dots, \dim X_\ar$.}
\vskip 0.20cm
To understand this general cohomology theory in down-to-earth terms, in section two, we develop a much more 
refined cohomology theory for Weil divisors $D$ over arithmetic surfaces $X$. This, in addition, is based on a basic 
theory for canonical ind-pro topologies over arithmetic adelic spaces. Recall that, by 
definition,
$$
\A_X^\ar:=\A_X^\ar(\O_X)\simeq\lim_{\raw D}\lim_{\substack{\law E\\ E\leq D}}\A_{X,12}^\ar(D)\big/\A_{X,12}^\ar(E).
$$
Here $\A_{X,12}^\ar(D)$ is one of the level two subspaces of $\A_X^\ar$ introduced in \S2.3.1. Moreover, for divisors 
$E\leq D$, 
$\A_{X,12}^\ar(D)\big/\A_{X,12}^\ar(E)$ are locally compact. Thus, using first projective 
then inductive limits, we obtain a canonical ind-pro topology on $\A_X^\ar$. In particular, we have the following 
natural generalization of topological theory for one dimensional adeles (see e.g., [Iw], [T]) to dimension two. 

\vskip 0.20cm
\noindent
{\bf Theorem II.}
{\it Let $X$ be an arithmetic surface. With respect to the canonical ind-pro topology on $\A_X^\ar$, we have

\noindent
(1) $\A_X^\ar$ is a Hausdorff, complete, and compact oriented topological group;

\noindent
(2) $\A_X^\ar$ is self-dual. That is, as topological groups,}
$$
\widehat{\ \A_X^\ar\ }\ \simeq\ \A_X^\ar.
$$\\[-2.6em]

With these basic topological structures exposed, next we apply them to our cohomology groups.  
For this, we first recall an arithmetic residue theory in \S2.1, by adopting a very precise approach of Morrow [M1,2], a 
special case of a general theory on residues of Grothendieck  (see e.g., [L], [B] and [Y]). Then, we introduce a 
global pairing $\langle\cdot,\cdot\rangle_\omega$ in \S2.2 on the arithmetic adelic space $\A_X^\ar$, and prove the 
following

\vskip 0.20cm
\noindent
{\bf Proposition A.}
{\it Let $X$ be an arithmetic surface and $\omega$ be a non-zero rational differential on $X$. Then, the natural residue
pairing $\dis{\langle\cdot,\cdot\rangle_\omega:\A_X^\ar\times\A_X^\ar\to\mathbb S^1}$ is non-degenerate.} 

Moreover, we construct the so-called level two adelic subspaces $\A_{X,01}^\ar, \A_{X,02}^\ar$ and 
$\A_{X,12}^\ar(D)$ of $\A_X^\ar$ in \S2.3.1. Accordingly, we calculate their perpendicular subspaces with respect to our 
global residue pairing.

\vskip 0.20cm
\noindent
{\bf Proposition B.}
{\it Let $X$ be an arithmetic surface, $D$ be a Weil divisor and $\omega$ be a non-zero rational differential on $X$. 
Denote by $(\omega)$ the canonical divisor on $X$ associated to $\omega$. Then we have

\noindent
(1) Level two subspaces $\A_{X,01}^\ar, \A_{X,02}^\ar$ and $\A_{X,12}^\ar(D)$ are closed in $\A_X^\ar$;

\noindent
(2) With respect to the residue pairing} $\langle\cdot,\cdot\rangle_\omega$,
$$
\big(A_{X,01}^{\ar}\big)^\perp=A_{X,01}^{\ar},\quad \big(A_{X,02}^{\ar}\big)^\perp=A_{X,02}^{\ar},
\quad\mathrm{and}\quad \big(A_{X,12}^{\ar}(D)\big)^\perp=A_{X,12}^{\ar}((\omega)-D).
$$

Our lengthy proof for (2) is based on the residue formulas for horizontal and vertical curves on arithmetic surfaces 
established in [M2]. Moreover, as one can find from the proof of this theorem in \S2.3.2, all the level two adelic 
subspaces $\A_{X,01}^\ar, \A_{X,02}^\ar$ and $\A_{X,12}^\ar(D)$ are characterized by these perpendicular properties 
as well. As for (1), our proof in \S 3.1.3 uses a topological notion of completeness in an essential way.
\vskip 0.10cm
With the help of these level two subspaces, now we are ready to write down the adelic complex  of \S1.2.3 and 
hence its cohomology groups $H_\ar^i(X,\O_X(D))$ associated to the  line bundle $\O_X(D)$ on an arithmetic surface 
$X$ very explicitly. Indeed, according to \S1.2.3, or more directly, following [P], we arrive at the following central
\vskip 0.20cm
\noindent
{\bf Definition.} 
{\it Let $X$ be an arithmetic surface and $D$ be a Weil divisor on $X$. We define arithmetic cohomology groups 
$H_\ar^i(X,\O_X(D))$ for the line bundle $\O_X(D)$ on $X$, $i=0,1,2$,  by}
$$
\begin{aligned}
H^0_{\ar}(X,\O_X(D)):=&\A_{X, 01}^{\ar}\cap \A_{X, 02}^{\ar}\cap \A_{X, 12}^{\ar}(D);\\[0.2em]
H^1_{\ar}(X,\O_X(D))&\\
:=\Big(\big(\A_{X, 01}^{\ar}+\A_{X, 02}^{\ar}&\big)\cap \A_{X, 12}^{\ar}(D)\Big)
\Big/\Big(\A_{X, 01}^{\ar}\cap \A_{X, 12}^{\ar}(D)+\A_{X, 02}^{\ar}\cap \A_{X, 12}^{\ar}(D)\Big);\\[0.2em]
H^2_{\ar}(X,\O_X(D)):=&\A_{X, 012}^{\ar}\Big/\Big(\A_{X, 01}^{\ar}+ \A_{X, 02}^{\ar}+\A_{X, 12}^{\ar}(D)\Big).
\end{aligned}
$$

Similar to usual cohomology theory, these cohomology groups admit a natural inductive structure. For details, please 
refer to Propositions 17, 18 in \S 2.4.2. Moreover, induced from the canonical ind-pro topology on $\A_X^\ar$, we obtain 
natural topological structures for our cohomology groups, since, from Proposition B(1) above, the subspaces 
$\A_{X,01}^\ar, \A_{X,02}^\ar$ and $\A_{X,12}^\ar(D)$ are all closed. Consequently, as one of the main results of this 
paper, with the use of Theorem II above, in \S 3.2.4, we are able to establish the following
\vskip 0.20cm
\noindent
{\bf Theorem III.}
{\it Let $X$ be an arithmetic surface with a canonical divisor $K_X$ and $D$ be a Weil divisor on $X$. Then, as 
topological groups,} 
$$
\widehat {H_\ar^i(X,\O_X(D))}\simeq H_\ar^{2-i}(X,\O_X(K_X-D))\qquad i=0,1,2.
$$
 
Our theory is natural and proves to be very useful. For example, as recalled in \S1.2.4, in [W], based on Tate's 
thesis ([T]), for a metrized bundle $\overline E$ on an arithmetic curve $\Spec\, \O_F$, we are able to prove a refined 
arithmetic duality: 
$$
h^1_{\ar}\big(X,\, \overline E\big)\,=\,h^0_{\ar}\big(X,\ \overline {K_X}\otimes \overline E^\vee\big),
$$ 
and establish \lq the' arithmetic Riemann-Roch theorem:
$$
h^0_{\ar}\big(X,\, \overline E\big)-h^1_{\ar}\big(X,\, \overline E\big)
=\deg_\ar(\overline E)-\frac{\mathrm{rank}E}{2}\log|\Delta_F|,
$$ 
(where $h^i$ denotes the arithmetic count of $H^i_{\ar}$,)  and obtain an effective version of ampleness and 
vanishing theorem. All this plays an essential role in our studies of non-abelian zeta functions for number fields.

\section{Arithmetic Adelic Complexes}

\subsection{Parshin-Beilinson's Theory}

For later use, we here recall some basic constructions of adelic cohomology theory for Noetherian schemes of 
Parshin-Beilinson ([P1,2], [B], see also [H]).

\subsubsection{Local fields for reduced flags}
Let $F$ be a number field with $\mathcal O_F$ the ring of integers, and $\pi: X\to \Spec\, \O_F$ be an integral  
arithmetic variety. By a flag $\delta=(p_0,p_1,\dots, p_n)$ on $X$, we mean a chain of integral subschemes $p_i$ 
satisfying $\displaystyle{p_{i+1}\in\overline{\{p_i\}}=:X_i}$; and we call $\delta$ reduced if $\dim p_i=n-i$ for each $i$. 
For a reduced $\delta$, with respect to each affine open neighborhood $U=\Spec\, B$ of the closed point $p_n$, we 
obtain a chain, denoted also by $\delta$ with an abuse of notation, of prime divisors on $U$. Consequently, through
processes of localizations  and completions, we can associate to $\delta$ a ring 
$$
\displaystyle{K_\delta:=C_{p_0}S_{p_0}^{-1}\dots C_{p_n}S_{p_n}^{-1}B}.
$$
Here, as usual, for a ring $R$, an $R$-module $M$ and a prime ideal $p$ of $R$, we write $S_p^{-1}M$ for the 
localization of $M$ at $S_p=R\backslash p$, and $C_pM=\lim_{\leftarrow_{n\in\mathbb N}}M/p^nM$ its $p$-adic 
completion.

The ring $K_\delta$ is independent of the choices of $B$. Indeed, following [P2], we can introduce inductively
schemes $X_{i,\alpha_i}'$ as in the following diagram
$$
\begin{matrix}X_0&\supset&X_1&\supset&X_2&\supset&\cdots\\
\uparrow&&\uparrow&&&&\\
X_0'&\supset&X_{1,\alpha_1}&\supset&\uparrow&&\\
&&\uparrow&&&&\\
&&X_{1,\alpha_1}'&\supset&X_{2,\alpha_2}&&\\
&&&&\uparrow&&\\
&&&&\vdots&&\\\end{matrix}
$$
where $X'$ denotes the normalization of a scheme $X$, and $X_{i,\alpha_i}$ denotes an integral subscheme in 
$X_{i-1,\alpha_{i-1}}'$ which dominates $X_i$. In particular, 

\noindent
(i) $X_{1,\alpha_1}$, being an integral  subvariety of the normal scheme $X_0'$, defines a discrete valuation of the 
field of rational functions on $X_0$, whose residue field coincides with the field of rational functions on the normal 
scheme $X_{1,\alpha_1}'$. 

\noindent
(ii) More generally, for a fixed  $i,\, 1\leq i\leq n$, $X_{i,\alpha_i}$, being an integral  subvariety of the normal scheme
$X_{i-1,\alpha_{i-1}}'$, defines a discrete valuation of the field of rational functions on $X_{i-1,\alpha_{i-1}}'$, whose
residue field coincides with the field of rational functions on the normal scheme $X_{i,\alpha_i}'$.

\noindent
Accordingly, for each collection $(\alpha_1,\alpha_2,\dots,\alpha_n)$ of indices, the chain of field of rational functions 
$K_0,K_{1,\alpha_1},\dots, K_{n,\alpha_n}$ defines an $n$-dimensional local field $K_{(\alpha_1,\dots,\alpha_n)}$
and hence an Artin ring
$$
\mathbb K_\delta:=\oplus_{(\alpha_1,\dots,\alpha_n)\in\Lambda_\delta}K_{(\alpha_1,\dots,\alpha_n)}.
$$

\begin{thm} ([P1,2], [Y]) Let $\delta=(p_0,p_1,\dots,p_n)$ be a reduced flag on $X$, and 
$K_{(\alpha_1,\dots,\alpha_n)}$ be the $n$-dimensional local field associated to the collection of indices
$(\alpha_1,\alpha_2,\dots,\alpha_n)$ above. Then,  we have

\noindent
(1) The $n$-dimensional local field $K_{(\alpha_1,\dots,\alpha_n)}$ is, up to finite extension, isomorphic to 
$$
{either}\quad F_v'((t_{n-1}))\cdots((t_1)),\quad {or}\quad 
F_v'\{\{t_n\}\}\cdots\{t_{m+2}\}\}((t_m))\dots ((t_1))
$$ 
where $F_v'$ denotes a certain finite extension of some $v$-adic non-Archimedean local field $F_v$;

\noindent
(2)  The ring $K_\delta$ is isomorphic to $\mathbb K_\delta$. In particular, it is independent of the choices of $U$.
\end{thm}

For example, if $X$ is an arithmetic surface, and $p_1$ is  a vertical curve, then, up to finite extension,
$K_\delta=F_{\pi(p_2)}\{\{u\}\}$,\footnote{Definition of $F_{\pi(p_2)}\{\{u\}\}$ will be recalled in \S2.1.1.}
where $u$ denotes a local parameter of the curve $p_1$ at the point $p_2$, and $F_{\pi(p_2)}$ denotes the 
$\pi(p_2)$-adic number field  associated to the closed point $\pi(p_2)$ on $\Spec\,\O_F$; on the other hand, if $p_1$ is a 
horizontal curve, then $K_\delta=L((t))$, where $t$ is a local parameter of $p_1$ at $p_2$, and $L/F$ is a finite field 
extension. Indeed, $p_1$ corresponds to an algebraic point on the generic fiber $X_F$ of $\pi$, and $L$ is simply 
the corresponding defining field.

\subsubsection{Adelic cohomology theory}

Let $X$ be a Noetherian scheme, and let $P(X)$ be the set of (integral) points of $X$ (in the scheme theoretic sense). 
For $p,\,q\in P(X)$,  if $q\in\overline{\{p\}}$, we write $p\geq q$. Let $S(X)$ be the simplicial set induced by 
$(P(X),\geq)$, i.e., the set of $m$-simplices of $S(X)$ is defined by
$$
S(X)_m:=\big\{(p_0,\dots,p_m)\,|\ p_i\in P(X),\ p_i\geq p_{i+1}\big\},
$$ 
the natural boundary maps $\delta_i^n$ are defined by deleting the $i$-th point, and the degeneracy maps 
$\sigma_i^n$ are defined by duplicating the $i$-th point:
$$
\begin{aligned}
\delta_i^m:S(X)_m\to S(X)_{m-1},&\ \  (p_{0}, \dots, p_{i}, \dots, p_{m})\mapsto(p_0,\dots,\check{p_i},\dots,p_m),\\[0.60em]
\sigma_i^m:S(X)_m\to S(X)_{m+1},&\ \ (p_{0}, \dots, p_{i}, \dots, p_{m})\mapsto (p_0,\dots,p_i, p_i,\dots,p_m).\\
\end{aligned}
$$
Denote also by $S(X)^{\mathrm{red}}_m$ the subset of $S(X)_m$ consisting of all non-degenerate simplexes, i.e.,
$$
S(X)_m^{\mathrm{red}}=\big\{(p_0,\dots,p_m)\in S(X)_m\,\big|\ \dim p_i\not =\dim p_j\ \, \forall i\not= j\,\big\}.
$$
For $p\in P(X)$ and $M$ an $\mathcal O_p$-module, set $[M]_p:=(i_p)_*M$, where 
$i_p:\Spec(\mathcal O_p)\hookrightarrow X$ denotes the natural induced morphism. Moreover, for 
$K\subset S(X)_m$ and a point $p\in P(X)$, introduce $_pK\subset S(X)_{m-1}$ by
$$
_pK:=\big\{(p_1,\dots,p_m)\in S(X)_{m-1}\,|\ (p,p_1,\dots,p_m)\in K\big\}.
$$
Then, we have the following

\begin{prop} ([P1,2], [B], see also [H, Prop 2.1.1]) There exists  a unique system of functors 
$\displaystyle{\{\mathbb A(K,*)\}_{_{K\subset S(X)}}}$ from the category of quasi-coherent sheaves on $X$ to the 
category of abelian groups, such that

\noindent
(i) $\mathbb A(K,\cdot)$ commutes with direct limits.

\noindent
(ii) For a coherent sheaf $\mathcal F$ on $X$,
$$
\mathbb A(K,\FF)=
\begin{cases}{\displaystyle\prod_{\ p\in K\ }\, \ {\lim_{\longleftarrow}}_l}\,\FF_p\big/\frak m_p^l\FF_p,&m=0,\\[1.80em]
{\displaystyle\prod_{p\in P(X)}{{\lim_{\longleftarrow}}_l}}\,\A\big(\,_pK,[\FF_p\big/\frak m_p^l\FF_p]_p\big),&m>0.
\end{cases}
$$
Here  $\frak m_p$ denotes the prime ideal associated to $p$. 
\end{prop}

Consequently, for any quasi-coherent sheaf $\FF$ on $X$, there exist well-defined adelic spaces 
$$
\A_X^m(\FF):=\A\big(S(X)_m^{\mathrm{red}},\FF\big).
$$
Clearly, if we introduce $K_{i_0,\dots, i_m}
=\big\{(p_0,\dots, p_m)\in S(X)_m\,|\,\mathrm{codim}\big(\overline{\{p_r\}}\big)=i_r\ \forall\, 0\leq r\leq m\big\}$, and 
define $\A_{X; i_0,\dots, i_m}(\FF):=\A_X\big(K_{i_0,\dots, i_m},\FF\big)$, then 
$$
\A_X^m(\FF)=\bigoplus_{0\leq i_0<\dots<i_m\leq\dim X}\A_{X; i_0,\dots, i_m}(\FF).
$$
Moreover, since $\A(K,\FF)\subset\prod_{(p_0,\dots,p_m)\in K}\A\big((p_0,\dots,p_m),\FF\big)$, we sometimes write 
an element $f$ of $\A(K,\FF)$ as $f=(f_{p_0,\dots,p_m})$ or $f=(f_{X_0,\dots,X_m})$, where $X_i=\overline{\{p_i\}}$
and $f_{p_0,\dots,p_m}=f_{X_0,\dots,X_m}\in \A\big((p_0,\dots,p_m),\FF\big)$.

To get an adelic complex associated to $X$, we next introduce boundary maps 
${d^m: \A_X^{m-1}(\FF)\to\A_X^{m}(\FF)}$ as in [H, Def 2.2.2]. For $K\subset S(X)_m$ and $L\subset S(X)_{m-1}$ 
such that $\delta_i^mK\subset L$ for a certain $i$, we define a boundary map 
$${d_i^m(K,L,\FF):\ \A(L,\FF)\longrightarrow \A(K,\FF)}$$ as follows. 

\noindent
(a) For coherent sheaves $\FF$, 

\noindent
(i) When $i=0$, for $p\in P(X)$, induced from the morphism $\FF\to [\FF_p/\frak m_p^l\FF_p]_p$ and the inclusion
$_pK\subset L$, we have
the morphisms $\A(L,\FF)\to \A\big(L,[\FF_p/\frak m_p^l\FF_p]_p\big)$ and 
$\A\big(L,[\FF_p/\frak m_p^l\FF_p]_p\big)\to  \A\big(\,_pK,[\FF_p/\frak m_p^l\FF_p]_p\big)$. Their compositions form a 
projective system $\varphi_p^l:\A(L,\FF)\to \A\big(\,_pK,[\FF_p/\frak m_p^l\FF_p]_p\big).$ Accordingly, we set 
$d_0^m(K,L,\FF):=\prod_{p\in P(X)}\lim_{\longleftarrow l}\varphi_p^l;$

\noindent
(ii) When $i=m=1$, we obtain a projective system induced from the standard morphisms 
$\pi_p^l:\Gamma\big(X,[\FF_p/\frak m_p^l\FF_p]_p\big)\to\A\big(\,_pK,[\FF_p/\frak m_p^l\FF_p]_p\big)$. Accordingly, 
we set $d_1^1(K,L,\FF):=\prod_{p\in P(X)}\lim_{\longleftarrow l}\pi_p^l;$

\noindent
(iii) When $i>0, m>0$, we use an induction on $(i,m)$. That is to say, we set 
$d_i^m(K,L,\FF):=\prod_{p\in P(X)}\lim_{\longleftarrow l}d_{i-1}^{m-1}\big(\,_pK,\,_pL,[\FF_p/\frak m_p^l\FF_p]_p\big).$

\noindent
(b) For quasi-coherent sheaves $\FF$, first we write $\FF$ as an inductive limit of coherent sheaves, then we use (a) 
to get boundary maps for the later, finally we use the fact that in the definition of (a), all constructions commute with 
inductive limits. One checks (see e.g. [H]) that the resulting boundary map is well-defined.

With this, set
$$
d^m:=\sum_{i=0}^m(-1)^i\,d_i^m\big(S(X)_m^{\red},S(X)_{m-1}^{\red};\FF\big).
$$ 
Then we have the following

\begin{thm}  ([P1,2], [B], see also [H, Thm 4.2.3])  Let $X$ be a Noetherian scheme. Then, for  any quasi-coherent 
sheaf $\FF$ over $X$, we have

\noindent
(1)  $\big(\A_X^*(\FF),d^*\big)$ forms a cohomological complex of abelian groups;

\noindent
(2) Cohomology groups of  the complex $\big(\A_X^*(\FF),d^*\big)$ coincide with Grothendieck's sheaf 
theoretic cohomology groups $H^i\big(X,\FF\big)$. That is to say,  we have, for all $i$, 
$$
H^i\big(\A_X^*(\FF),\,d^*\big)\simeq H^i\big(X,\FF\big).
$$
\end{thm}

\subsubsection{Examples}
{\bf (A) Algebraic Curves}

Let $X$ be an integral regular projective curve defined over a field $k$. Denote its generic point by $\eta$ and its field 
of rational functions by $k(X)$. For a divisor $D$ on $X$, let $\O_X(D)$ be the associated invertible sheaf. Then,  from 
definition, the associated adelic spaces can be calculated as follows:
$$
\begin{aligned}
\A_{X;0}(\O_X(D))&=\A\big(\{\eta\},\O_X(D)\big)\\
=&\lim_{\longleftarrow_l}\O_X(D)_\eta\big/\frak m_\eta^l\O_X(D)_\eta
=\lim_{\longleftarrow_l}k(X)/\{0\}=k(X),
\end{aligned}
$$\\[-2.0em]
$$\begin{aligned}
&\A_{X;1}(\O_X(D))=\A_X\big(\{p\}\,|\, p\in X:\, \mathrm{closed\ point}\},\O_X(D)\big)\\
=&\prod_{p\in X}\lim_{\longleftarrow l}\O_X(D)_p\big/\frak m_p^l\O_X(D)_p
=\prod_{p\in X}\lim_{\longleftarrow l}\frak m_p^{-\ord_p(D)}\big/\frak m_p^{-\ord_p(D)+l}\\
=&\prod_{p\in X}\frak m_p^{-\ord_p(D)}=\Big\{(a_p)\in\prod_{p\in X}k(X)_p\,\big|\,\ord_p(a_p)+\ord_p(D)\geq 0\Big\},
\end{aligned}
$$
and
$$
\begin{aligned}
&\A_{X;01}(\O_X(D))=\A_X\big(\{\eta,p\,|\, p\in X:\, \mathrm{closed\ point}\},\O_X(D)\big)\\
=&\lim_{\longleftarrow l}\A\big(\{p\,|\, p\in X:\, \mathrm{closed\ point}\},
[\O_X(D)_\eta\big/\frak m_\eta^l\O_X(D)_\eta]_\eta\big)\\
=&\A\big(\{p\,|\, p\in X:\, \mathrm{closed\ point}\},\widetilde{k(X)}\big)\\
=&\A\big(\{p\,|\, p\in X:\, \mathrm{closed\ point}\},\lim_{\ \ \longrightarrow E}\O_X(E)\big)\\
=&\lim_{\ \, \longrightarrow E}\A_{X;1}(\O_X(E))=\bigcup_{E}\A_{X;1}(\O_X(E))\\
=&\Big\{(a_p)\in\prod_{p\in X}k(X)_p\,\big|\,a_p\in \O_p\ \forall' p\Big\}.
\end{aligned}
$$

\noindent
{\it Remark.} 
To calculate $\A_{X;01}(\O_X(D))$, when dealing with the constant sheaf $[k(X)]_\eta$, we cannot use 
Proposition 2(ii) directly, since $[k(X)]_\eta$ is not coherent. Instead, above, we first expressed it as an inductive limit 
of coherent sheaves $\O_X(E)$ associated to divisors $E$, then get the result from the inductive limit of adelic 
spaces for $\O_X(E)$'s. Indeed, if we had used Proposition 2 (ii) directly, then we would have obtained simply 
$\A\big(\{p\}\,|\, p\in X:\, \mathrm{closed\ point}\}, [k(X)]_\eta\big)=\{0\}$, a wrong claim.

Clearly, $\A_{X;01}(\O_X(D))$ is independent of $D$. We will write it as $\A_{X;01}$, or simply $\A_X$.
Consequently,  the associated adelic complex\\[-0.8em] 
$$
0\longrightarrow \A_{X;0}(\O_X(D))\oplus\A_{X;1}(\O_X(D))\buildrel d^1\over\longrightarrow\A_{X;01}(\O_X(D))
\longrightarrow 0
$$\\[-1.0em]
is given by 
$$ 
0\longrightarrow k(X)\oplus\A_{X;1}(\O_X(D))\buildrel d^1\over\longrightarrow\A_X\longrightarrow 0
$$ 
where $d^1:(a_0,a_1)\mapsto a_1-a_0$. Therefore,
$$
\begin{aligned}
H^0\big(\A_X(\O_X(D))\big)=&\,k(X)\,\cap\, \A_{X:1}(\O_X(D)),\\[-0.1em]
H^1\big(\A_X(\O_X(D))\big)=&\,\A_X\big/\big(k(X)+\A_{X;1}(\O_X(D))\big).
\end{aligned}
$$\\[-2.0em]

Note that $\A_X(\O_X(D))$ is simply $\A_X(D)$ of\, [S, Ch. 2], or better, [Iw, \S4]. We have proved the following

\begin{prop} (See e.g., [S, Ch. 2], [Iw,\,\S4]) For a divisor $D$ over an integral regular projective curve defined over a 
field $k$, we have\\[-0.8em] 
$$
H^0\big(\A_X(\O_X(D))\big)=H^0\big(X,\O_X(D)\big),\qquad H^1\big(\A_X(\O_X(D))\big)=H^1\big(X,\O_X(D)\big).
$$
\end{prop}
\vskip 0.20cm
\noindent
{\bf (B) $\mathbb A_{01}$ For Surfaces} 

Let now $X$ be a two dimensional integral, regular, Noetherian scheme, viewed as the finite part of an arithmetic surface. To calculate $\mathbb A_{X;01}$, similarly, by Proposition 2,
we have
$$
\begin{aligned}
&\A_{X;01}=\A_{X;01}(\O_X)=\A_X\big(\{\eta,C\,|\, C\subset X:\, \mathrm{irreducible\ curve}\},\O_X\big)\\
=&\lim_{\longleftarrow l}\A\big(\{C\,|\, C\subset X:\, \mathrm{irreducible\ curve}\},
[\O_{X,\eta}\big/\frak m_\eta^l\O_{X,\eta}]_\eta\big)\\
=&\A\big(\{C\,|\, C\subset X:\, \mathrm{irreducible\ curve}\},\widetilde{k(X)}\big)\\
=&\A\big(\{C\,|\, C\subset X:\, \mathrm{irreducible\ curve}\},\lim_{\ \ \longrightarrow E}\O_X(E)\big)\\
=&\lim_{\ \ \longrightarrow E}\A\big(\{C\,|\, C\subset X:\, \mathrm{irreducible\ curve}\},\O_X(E)\big)\\
=&\lim_{\ \ \longrightarrow E}\prod_{C:\, C\subset X\, \mathrm{irred\,curve}}
\lim_{\longleftarrow l}\O_X(E)_C\big/\frak m_C^l\O_X(E)_C\\
=&\lim_{\ \ \longrightarrow E}\prod_{C:\, C\subset X\, \mathrm{irred\,curve}}
\frak m_C^{-\nu_C(E)}\,\widehat{\,\mathcal O\,}_{\!\!X,C}.
\end{aligned}
$$
In particular, components of an element of $\A_{X;01}$ are independent of closed points. Thus, with diagonal embedding $\A_{X;01}\hookrightarrow \A_{X;012}$, we may write it as $(f_{C})_C$ instead of $(f_{C,x})_{C,x}$.

\subsection{Arithmetic Cohomology Groups}

Let $F$ be a number field with $\O_F$ the ring of integers. Denote by $S_{\fin}$, resp. $S_\infty$, the collection of 
finite, resp. infinite, places of $F$. Write $S=S_{\fin}\cup S_\infty$. Let $\pi:X\to\Spec\,\O_F$ be an integral arithmetic 
variety of pure dimension $n+1$. That is, an integral Noetherian scheme $X$, a flat and proper morphism $\pi$ with generic 
fiber $X_F$ a projective variety of dimension $n$ over $F$. For each $v\in S$, we write $F_v$ the $v$-completion of 
$F$, and for each $\sigma\in S_\infty$, we write $X_v:=X\times_{\O_F}\Spec F_v$ and write 
$\varphi_\sigma:X_\sigma\to X_F$ for the map induced from the natural embedding $F\hookrightarrow F_\sigma$.
In particular, an arithmetic variety $X$ consists of two parts, the finite one, which we also denote by $X$, and an 
infinite one, which we denote by $X_\infty$. These two parts are closely interconnected. 

\subsubsection{Adelic rings for arithmetic surfaces} 

The part of our theory on arithmetic adelic complexes for finite places now becomes very simple. Indeed, our 
arithmetic variety $X$ is assumed to be Noetherian, so we can apply the theory recalled in \S1.1 directly. In particular, 
for a quasi-coherent sheaf $\FF$ on $X$, we have well-defined adelic spaces\\[-0.8em]  
$$
\A^{\fin}_{X;\,i_0,\dots,i_m}(\FF):=\A_X(K_{i_0,\dots,i_m},\FF).
$$

So to define $\A^{\ar}_{X;\,i_0,\dots,i_m}(\FF)$, we need to understand what happens on $X_\infty$. For this purpose,
we next recall Osipov-Parshin's construction of arithmetic adelic ring $\A_X^\ar$ for an arithmetic surface $X$.

\begin{defn} ([OP]) {\bf Arithmetic adelic ring of an arithmetic surface} Let $\pi:X\to\Spec\,\O_F$ be an arithmetic 
surface, i.e., a 2-dimensional arithmetic variety, with generic fiber $X_F$.

\noindent
(i) Finite adelic ring: From the Parshin-Beilinson theory for the Noetherian scheme $X$, we  define
$$
\A^{\fin}_X:=\A_{X;012}(\O_X)=\lim_{\substack{\longrightarrow\\ D_1}}
\lim_{\substack{\longleftarrow\\ D_2:\,D_2\leq D_1}}\A_{X;12}(D_1)\big/\A_{X;12}(D_1).
$$ 
Here $D_*$'s are divisors on $X$ and $\A_{X;12}(D):=\A_{X;12}(\O_X(D_*))$ for $*=1,2$;

\noindent
(ii) $\infty$-adelic ring: Associated to the regular integral curve $X_F$ over $F$, we obtain the adelic ring 
$$
\A_{X_F}:=\A_{X_F;01}(\O_{X_F})=\lim_{\substack{\longrightarrow\\ {D}_1}}
\lim_{\substack{\longleftarrow\\ {D}_2:\,{D}_2\leq {D}_1}}\A_{X_F;1}({D}_1)\big/\A_{X_F;1}({D}_1).
$$ 
Here ${D}_*$'s are divisors on $X_F$ and $\A_{X_F;1}({D}):=\A_{X_F;1}(\O_{X_F}({D}_*))$ for $*=1,2$.

By definition,
$$
\A_X^\infty:=\A_{X_F}\,\widehat{\,\otimes\,} _{\mathbb Q}\,\mathbb R
:=\lim_{\substack{\longrightarrow\\ {D}_1}}\lim_{\substack{\longleftarrow\\ {D}_2:{D}_2\leq {D}_1}}
\Big(\big(\A_{X_F;1}({D}_1)\big/\A_{X_F;1}({D}_1)\big)\,\otimes _{\mathbb Q}\,\mathbb R\Big).
$$

\noindent
(iii) Arithmetic adelic ring: The arithmetic adelic ring of an arithmetic surface $X$ is defined by
$$
\A_{X}^{\ar}:=\A_{X;012}^\ar:=\A_X^{\fin}\bigoplus\A_X^\infty.
$$
\end{defn}

The essential point here is, for divisors $D_i,\, i=1,2,$ over the curve $X_F$, when $D_2\leq D_1$, the quotient 
$\A_{X;1}({D}_1)\big/\A_{X;1}({D}_1)$ is a  finite dimensional $F$- and hence $\mathbb Q$-vector space. 

To help the reader understand this formal definition in concrete terms, we add following examples.

\begin{ex} 
On $X=\mathbb P^1_{\mathbb Z}$
\end{ex} 
\noindent
We have $X_{\mathbb Q}=\mathbb P^1_{\mathbb Q}$ and
$\mathbb Q\big(\mathbb P^1_{\mathbb Q}\big)=\mathbb Q(t)$. Easily, 
$$
\mathbb Q((t))\otimes_{\mathbb Q} {\mathbb R}\not=\mathbb R((t)).
$$
However, since $\displaystyle{\mathbb Q((t))=\lim_{\substack{\longrightarrow\\ n}}
\lim_{\substack{\longleftarrow\\ m:\,m\leq n}}t^{-n}\mathbb Q[[t]]\,\big/\, t^{-m}\mathbb Q[[t]]}$ and the 
$\mathbb Q$-vector spaces $\displaystyle{t^{-n}\mathbb Q[[t]]\,\big/\, t^{-m}\mathbb Q[[t]]}$ are finite dimensional,
we have 
$$
\begin{aligned}
\mathbb Q((t))\,\widehat{\,\otimes\,}_{\mathbb Q}\,\mathbb R
=&\lim_{\substack{\longrightarrow\\ n}}\lim_{\substack{\longleftarrow\\ m:\,m\leq n}}
\big(t^{-n}\mathbb Q[[t]]\big/ t^{-m}\mathbb Q[[t]]\big){\,\otimes\,}_{\mathbb Q}\,\mathbb R\\ 
=&\lim_{\substack{\longrightarrow\\ n}}\lim_{\substack{\longleftarrow\\ m:\,m\leq n}}
\big(t^{-n}\mathbb R[[t]]\big/ t^{-m}\mathbb R[[t]]\big)=\mathbb R((t)).
\end{aligned}
$$

\begin{ex} 
Over an arithmetic surface $X$
\end{ex}
\noindent
For a complete flag $(X, C,x)$ on $X$ (with $C$ an irreducible curve on 
$X$ and $x$ a close point on $C$), let $k(X)_{C,x}$ its associated local ring. By Theorem 1, $k(X)_{C,x}$ is a direct 
sum of two dimensional local fields. Denote by $\pi_C$ the local parameter defined by $C$ in $X$. Then 
$$
\begin{aligned}
\A_X^\fin=&\A_{X,012}={\prod_{x\in C}}'k(X)_{C,x}:={\prod_C}'\Big({\prod}_{x:x\in C}'k(X)_{C,x}\Big)\\
:=&\Big\{\big(\sum_{i_C=-\infty}^\infty h_C(a_{i_C})\pi_C^{i_C}\big)_C\in\prod_C\big(\prod_{x: x\in C}k(X)_{C,x}\big):\\
&\qquad a_{i_C}\in \A_{C,01},\ a_{i_C}=0\ (i_C\ll 0);\ \min\{i_C:a_{i_C}\not=0\}\geq 0\ (\forall' C)\Big\}
,\end{aligned}
$$ 
where $h_C$ is a lifting defined in [MZ], which we call the Madunts-Zhukov lifting.
For details, please see \S3.1.2. 

\subsubsection{Adelic spaces at infinity}

Now we are ready to treat adelic spaces at infinite places for general arithmetic varieties. Motivated by the discussion
above, we make the following

\begin{defn} Let $\pi:X\to\Spec\,\O_F$ be an arithmetic variety. Let $S(X_F)$ be the simplicial set associated to 
its generic fiber $X_F$ and $K\subset S(X_F)_m,\, m\geq 0$, a subset.

\noindent
(i) Let $\mathcal G$ be a coherent sheaf on $X_F$. We define the associated adelic spaces by  
$$
\A^\infty(K,\mathcal G):=
\begin{cases}
\displaystyle{\prod_{\ p\in K\ }{\lim_{\longleftarrow}}_l\big(\mathcal G_p/\frak m_p^l\mathcal G_p\otimes_{\mathbb Q} \R\big),}
&m=0\\[1.80em]
\displaystyle{\prod_{p\in P(X)}{\lim_{\longleftarrow}}_l\A^\infty\big(\,_pK,[\mathcal G_p/\frak m_p^l\mathcal G_p]_p\big),}&m>0.
\end{cases}
$$

\noindent
(ii) Let $\{\mathcal G_i\}_i$ be an inductive system of coherent sheaves on $X_F$ and 
$\displaystyle{\FF=\lim_{\longrightarrow i}\mathcal G_i}$. Then we define 
$$
\A^\infty(K,\FF):=\lim_{\longrightarrow i}\A^\infty(K,\mathcal G_i).
$$
\end{defn}
Clearly, the essential part of this definition is the one for $m=0$. Moreover, if 
$\displaystyle{\FF=\lim_{\longrightarrow i}\mathcal G_i'}$ is another inductive limit of coherent sheaves, we have 
$\displaystyle{\lim_{\longrightarrow i}\A^\infty(K,\mathcal G_i')\simeq \lim_{\longrightarrow i}\A^\infty(K,\mathcal G_i)}$,
by the universal property of inductive limits since $\mathcal G_p/\frak m_p^l\mathcal G_p$'s are 
$\mathbb Q$-vector spaces. Therefore, $\A^\infty(K,\FF)$ is well-defined for all quasi-coherent sheaves $\FF$ on $X_F$. 
Moreover, as a functor from the category of coherent sheaves on $X_F$ to that of $\mathbb Q$-vector spaces,  
$\A^\infty(K,*)$ is additive and exact. Hence, by [H, \S1.2], $\A^\infty(K,*)$ commutes with the direct  limits, even in 
general, for an inductive system $\{\FF_i\}_i$ of quasi-coherent sheaves 
$\ \displaystyle{\lim_{\longrightarrow i}}\,\A^\infty(K,\FF_i)\not=A^\infty(K,\displaystyle{\lim_{\longrightarrow i}}\,\FF_i).$

\subsubsection{Arithmetic adelic complexes}

As mentioned at the beginning of this section, for arithmetic varieties, the finite and infinite parts are closely interconnected. 
Therefore, when developing an arithmetic cohomology theory, we will treat them as an unify one using an uniformity condition.

Let  $X$ be an arithmetic variety with generic fiber $X_F$. For a point $P$ of $X_F$, denote its associated Zariski closure 
in $X$ by $E_{P}$. We call a flag $\delta=(\frak p_0,\frak p_1,\dots,\frak p_k)\in S(X)$ {\it horizontal}, if there exists a flag 
$\delta_F=(P_0,P_1,\dots, P_k)\in S(X_F)$ such that 
$(\overline{\frak p}_0,\overline{\frak p}_1,\dots,\overline{\frak p}_k)=(E_{P_0},E_{P_1},\dots, E_{P_k})$. Accordingly, for 
$K\subset S(X)$, we denote $K^{\mathrm{h}}$ the collection of all horizontal flags in $K$ and 
$K^{\mathrm{nh}}=K\smallsetminus K^{\mathrm{h}}$. Simply put, our uniformity condition is a constrain on adelic components 
associated to horizontal flags. 

Let $\FF$ be a quasi-coherent sheaf on $X$, denote its induced sheaf on the generic fiber $X_F$ by $\FF_F$. It is well-known 
that $\FF_F$ is quasi-coherent as well. Motivated by [W], we introduce the following

\begin{defn} Let $X$ be an arithmetic variety of dimension $n+1$ and $\FF$ a quasi-coherent  sheaf on $X$.
Fix an index tuple $(i_0,\dots, i_m)$ satisfying $i_0\leq \dots\leq i_m$.

\noindent
(i) The finite, resp. infinite, adelic space of type $(i_0,\dots, i_m)$ associated to $\FF$ is defined by
$$
\begin{aligned}
\A^{\fin}_{X;\,i_0,\dots, i_m}(\FF):=\A_X\big(K_{X;\,i_0,\dots, i_m},\FF\big)=
&\A_X^\fin\big(K^{\mathrm{nh}}_{X;\,i_0,\dots, i_m},\FF\big)\oplus
\A_X^\fin\big(K^{\mathrm{h}}_{X;\,i_0,\dots, i_m},\FF\big),\\[0.5em]
\quad{resp.}\qquad
\A^{\infty}_{X;\,i_0,\dots, i_m}(\FF):=&\A_X^\infty\big(K_{X_F;\,i_0,\dots, i_m},\FF_F\big).
\end{aligned}
$$
Here, for $Z\subset X$ or $X_F$, we set
$$K_{Z;i_0,\dots,i_m}:=\big\{(p_0,\dots, p_m)\in S(Z)_m\,\big|\,\mathrm{codim}_Z\,\overline{\{p_r\}}=i_r\ 
\forall\, 0\leq t\leq m\big\};$$

\noindent
(ii) The arithmetic adelic space of type $(i_0,\dots, i_m)$ associated to $\FF$ is defined by
$$
\A^\ar_{X;\,i_0,\dots,i_m}(\FF)\ \,=:
\begin{cases}\A^\fin_{X;\,i_0,\dots,i_m}(\FF)\,\bigoplus\, \A^\infty_{X;\,i_0,\dots,i_{m-1}}(\FF_F),& i_m=n+1;\\[0.5em]
\A_X^\fin\big(K^{\mathrm{nh}}_{X;\,i_0,\dots, i_m},\FF\big)\oplus 
\A_X^{\fin,\mathrm{inf}}\big(K^{\mathrm{h}}_{X;\,i_0,\dots, i_m},\FF\big),&i_m\not=n+1
\end{cases}
$$
where
$$
\A_X^{\fin,\mathrm{inf}}\big(K^{\mathrm{h}}_{X;\,i_0,\dots, i_m},\FF\big)\subset
\A_X^\fin\big(K^{\mathrm{h}}_{X;\,i_0,\dots, i_m},\FF\big)\,\bigoplus\, \A^\infty_{X;\,i_0,\dots, i_m}(\FF_F)
$$
consisting of adeles satisfying, for all flags 
$(\frak p_{i_0},\frak p_{i_1},\dots, \frak p_{i_m})\in K_{X_F;\,i_0,\dots, i_m}$,
$$
f_{E_{\frak p_{i_0}},E_{\frak p_{i_1}},\dots,E_{\frak p_{i_m}} }=f_{\frak p_{i_0},\frak p_{i_1},\dots, \frak p_{i_m}};
$$

\noindent
(iii) For $m\geq 0$, define the $m$-th reduced  arithmetic adelic space $\A^{\ar}_{X;\,m}(\FF)$ of $\FF$  by
$$
\A_{\ar,\,\red}^{m}(X,\FF):=\bigoplus_{\substack{(i_0,\dots,i_m)\\ 0\leq i_0<i_1<\dots<i_m\leq n+1}}\A^{\ar}_{X;\,i_0,\dots,i_m}(\FF).
$$
\end{defn}
 
\noindent
{\it Remarks}. (i) For any $\frak p\in P(X_F)$, $\O_{X,E_{\frak p}}=\O_{X_F,\frak p}$ and 
$k(X)_{E_{\frak p}}=k(X_F)_{\frak p}$. Consequently, for any $(\frak p_0,\dots,\frak p_m)\in S(X_F)_m$, we have a 
natural morphism 
$$
\A\big((E_{\frak p_0},\dots,E_{\frak p_m}),\FF\big)=\A\big((\frak p_0,\dots,\frak p_m),\FF_F\big).
$$
since $\FF$ is quasi-coherent. It is in this sense we use the relation
$f_{E_{\frak p_0},E_{\frak p_1},\dots,E_{\frak p_{m}} }=f_{\frak p_0,\frak p_1,\dots, \frak p_{m}}$ above. (In particular, if 
$\frak p_i$'s are vertical, there are no conditions on the corresponding components.) Clearly, this uniformity condition 
is an essential one, since it characters the natural interconnection between  finite  and infinite components of 
arithmetic adelic elements.

\noindent
(ii) In part (ii) of the definition, we need the space $\A^\infty_{X;\,\emptyset}(\FF_F)$. Here, to complete our definition,  for an 
arithmetic variety $X$, we view $\A^\infty_{X;\,\emptyset}(\FF_F)$ as the $(-1)$-level of the adelic complex for its generic fiber 
$X_F$. That is to say, we define it as follows. By [Y, p. 63], we have the (-1)-simplex $\underline{1}_{U}$ for open 
$U\subset X$. Set then $S(X_{F})_{-1}=\{\underbar{1}_{U}\mid U\subset X:{\rm open}\}$,  and, for $K\subset S(X_F)_{-1}$, let 
$$
\A^\infty_{X;\,\emptyset}(K, \FF_F):=
\begin{cases}
\FF_F(U_{K,F})\otimes_{\mathbb Q}\mathbb R,&\dim X\geq 2\\
\{\,s_F\in \FF_F(U_{K,F})\otimes \R\, |\, s\in\FF(U_K)\,\},&\dim X=1
\end{cases}
$$
where $U_K:=\cup_{\underline{1}_U\in K}U$ and $s_F$ denotes the section induced by $s$. The reason for separation of 
arithmetic curves with others in this latest definition is that arithmetic varieties are relative over arithmetic curves.

Moreover, from standard homotopy theory, if we introduce the boundary morphisms by 
$$
\begin{matrix}
{d^m_{i}:}&\bigoplus\ \A^{\ar}_{X;\,l_0,\dots,l_{m-1}}(\FF)&\longrightarrow&
\bigoplus\ \A^{\ar}_{X;\,k_0,\dots,k_{m}}(\FF)\\
&(a_{l_0,\dots,l_{m-1}})&\mapsto& (a_{k_0,\dots,\hat k_i,\dots,k_{m}});
\\\end{matrix}
$$ 
and $\ d_{m}=\sum_{i=0}^{m}(-1)^id^{m}_i:\  
\bigoplus\ \A^{\ar}_{X;\,k_0,\dots,l_{m-1}}(\FF)\ \ \longrightarrow \ \  \bigoplus\  \A^{\ar}_{X;\,k_0,\dots,k_{m}}(\FF),$
we have

\begin{prop}
$\Big(\A_{\ar,\,\red}^*(X,\FF),d^*\Big)$ defines a complex of abelian groups.
\end{prop}

All in all, we are now ready to introduce the following
\vskip 0.20cm
\noindent
{\bf Main Definition.} 
{\it Let $\pi: X\to\Spec\,\O_F$ be an arithmetic variety. Let $\FF$ be a quasi-coherent sheaf on 
$X$. Then we define the $i$-th adelic arithmetic cohomology groups of $\FF$ by
$$
H^i_{\ar}(X,\FF):=H^i(\A_{\ar,\,\red}^*(X,\FF),d^*\Big),
$$ 
the $i$-th cohomology group of the complex $(\A_{\ar,\,\red}^*(X,\FF),d^*\Big)$.}

\vskip 0.20cm
Consequently, we have the following

\begin{thm} If $X$ is an arithmetic variety of dimension $n+1$, then 
$$\hskip 2.0cm 
H^i_{\ar}(X,\FF)=0
\qquad\mathrm{unless}\ i=0,1,\dots, n+1.$$
\end{thm}

\noindent
{\it Proof.} Indeed, outside the range $0\leq i\leq n+1$, the complex consists of zero.

\subsubsection{Cohomology theory for arithmetic curves}

We here give an example of the above theory for arithmetic curves, which was previously developed in [W], based on Tate's  
thesis ([T]).
\vskip 0.10cm
Let $D=\sum_{i=1}^r n_i\frak p_i$ be a divisor on $X=\Spec\,\O_F$. Write $n_i=\ord_{\frak p_i}(D)$. For simplicity,
we use $\A_{\bullet}^*(D)$ instead of $\A_{\bullet}^*(\O_X(D))$. Then, by the same calculation as in \S1.1.3, we have
$\displaystyle{
\A^\fin_{X;01}(D)
=\big\{(a_p)\in{\prod}_{p\in X}F_p\,\big|\,a_p\in\O_p \ \forall' p\big\}.}
$
And, since $D_\eta=0$ is trivial, 
$
\A^\infty_{X;0}(D)
=\lim_{\leftarrow\, l}\O_{X_F,\eta}\big/\frak m_\eta^l\O_{X_F,\eta}\otimes_{\mathbb Q}\mathbb R
=\lim_{\leftarrow\, l}\big(F/\{0\}\big)\otimes_{\mathbb Q}\mathbb R
=F\otimes_{\mathbb Q}\mathbb R=\prod_{\sigma\in S_\infty}F_\sigma.
$
Therefore, 
$$
\A_{X;01}^\ar\big(\O_X(D)\big)
=\A^\fin_{01}(D)\oplus \A^\infty_0(D)
=\big\{ (a_{\frak p})\in{\prod}_{\frak p\in S}F_{\frak p}\,\big|\,a_{\frak p}\in \O_{\frak p}\ \forall' \frak p\in S_\fin\,\big\}.
$$
In particular, it coincides with the standard adelic ring $\A_F$ of $F$, hence is independent of $D$.
 
To understand $\A_0^\ar(D)$, we first calculate $\A_0^\fin(D)$. With the same calculation as in \S1.1.3 again, we have
$\A_{X;0}^\fin(D)=F$. Note that, from above, $\A^\infty_{X;0}(D)=F\otimes_{\mathbb Q}\mathbb R$. Thus, by 
definition, 
$$
\A_{X;0}^\ar(D)
=\big\{(a_v;a_\sigma)\in \A_{X;0}^\fin(D)\oplus \A^\infty_{X;0}(D)\,\big|\,
(a_v)=i_\fin(f), (a_\sigma)=i_\infty(f)\,\exists f\in F\big\}
$$
is then isomorphic to $F$, and hence also independent of $D$.
\vskip 0.15cm
From our definition, $\A_{X;1}^\ar(D)=\A^\fin_{X;1}(\O_X(D))\oplus\A_{X;\emptyset}^\infty(\O_X(D))$. To understand it, 
we first calculate $\A_{X;1}^\fin(D)$. With the same calculation as in \S1.1.3, we have
$\A_{X;1}^\fin(D)=\big\{(a_{\frak p})\in \A_{X;01}^\fin\,\big|\,\ord_{\frak p}(a_{\frak p})+\ord_{\frak p}(D)\geq 0\big\}.$
Then, by definition, we have
$\A_{X;\emptyset}^\infty(D)=\{s\in F \,\big|\,\ord_{\frak p}(s)+\ord_{\frak p}(D)\geq 0\big\},$
since $D=\sum_{\frak p}\ord_{\frak p}(D)\frak p$, and hence if $U=X-\{\frak p_1,\dots,\frak p_r\}$, $\O_X(U)$ is trivial.

In this way, we get the associated arithmetic adelic complex\\[-0.80em]
$$
0\longrightarrow \A^\ar_{X;0}(\O_X(D))\oplus\A^\ar_{X;1}(\O_X(D))\buildrel d^1\over\longrightarrow
\A^\ar_{X;01}(\O_X(D))\longrightarrow 0
$$
is given by:
$0\to F\oplus\A_{X;1}^\ar(\O_X(D))\buildrel d^1\over\longrightarrow\A_F\to 0, (a_0,a_1)\mapsto a_1-a_0.$ Therefore,
$$
\begin{aligned}
H_\ar^0\big(F,\O_X(D)\big)=&\,F\,\cap\, \A^\ar_{X;1}(\O_X(D)),\\[0.30em]
H_\ar^1\big(F,\O_X(D)\big)=&\,\A_F\big/\big(F+\A^\ar_{X;1}(\O_X(D))\big).
\end{aligned}
$$

In fact, a complete cohomology theory is developed for arithmetic curves $\Spec\,\O_F$ in [W]. 
For an $\O_F$-lattice $\Lambda$, i.e., a motorized locally free sheaf on $\Spec\,\O_F$, we introduce the associated topological cohomology groups $H^0_{\ar}(F, \Lambda)$ and 
$H^1_{\ar}(F,\Lambda)$, with $H^0_{\ar}(F, \Lambda)$ discrete and $H^1_{\ar}(F, \Lambda)$ compact. Consequently, using 
Fourier analysis for locally compact groups, we obtain their arithmetic counts $h^0_{\ar}(F, \Lambda)$ and 
$h^1_{\ar}(F,\Lambda)$. 

\begin{thm}
{\bf Cohomology Theory for Arithmetic Curves} ([W])\\[0.5em]
Let $F$ be a number field with $\O_F$ the ring of integers. Let $\omega_F$ be the Arakelov 
dualizing lattice of $\Spec\,\O_F$ and $\Delta_F$ be the discriminant of $F$. Then, for an $\O_F$-lattice $\Lambda$ 
of rank $n$ with its dual $\Lambda^\vee$, we have

\noindent
(1) (1.i) ({\bf Topological Duality}) As locally compact topological groups,
$$
\widehat{H^1_{\ar}(F,\omega_F\otimes \Lambda^\vee)}=H^0_{\ar}(F,\Lambda);
$$
\qquad\qquad Here $\,\widehat~\,$ denotes the Pontryagin dual.

(1.ii) ({\bf Arithmetic Duality})
$$
h^1_{\ar}(F,\omega_F\otimes \Lambda^\vee)=h^0_{\ar}(F,\Lambda);
$$

\noindent
(2) ({\bf Arithmetic Riemann-Roch Theorem}) 
$$
h^0_{\ar}(F,\Lambda)-h^1_{\ar}(F,\Lambda)=\mathrm{deg}_{\ar}(\Lambda)-\frac{n}{2}\log|\Delta_F|.
$$

\noindent
(3) ({\bf Ampleness, Positivity and Vanishing Theorem}) The following statements are equivalent:
 
(3.i) Rank one $\O_F$-lattice $A$ is arithmetic positive;

(3.ii) Rank one $\O_F$-lattice $A$ is arithmetic ample; and

(3.iii) For rank one $\O_F$-lattice $A$ and any $\O_F$-lattice $L$,
$$
\lim_{n\to\infty}h^1_\ar(F,A^n\otimes L)=0.
$$

\noindent
(4)  ({\bf Effective Vanishing Theorem}) Assume that $\Lambda$ is a semi-stable $\O_F$-lattice satisfying 
$\displaystyle{\deg_\ar(\Lambda)\leq-[F:\mathbb Q]\cdot\frac{n\log n}{2}}$,  then we have
$$
h^0_\ar(F,\Lambda)\leq \frac{\,3^{n\,[F:\mathbb Q]}}{1-\log 3/\pi}\cdot 
\exp\Big({-\pi[F:\mathbb Q]\cdot e^{-\frac{\deg_\ar(L)}{n}}}\Big).
$$ 
\end{thm}

For details, please refer to [W].

\eject
\section{Arithmetic Surfaces}

In the sequel, by an arithmetic surface, we mean a 2-dimensional regular integral Noetherian scheme $X$ together 
with a flat, proper morphism $\pi:X\to\Spec\,\O_F$. Here $\O_F$ denotes the ring of integers of a number field $F$. In 
particular, the generic fiber $X_F$ is a geometrically connected, regular, integral projective curve defined over $F$.

\subsection{Local Residue Pairings} 

Theory of residues for arithmetic surfaces, as a special case of Grothendieck's residue theory, can be realized using 
K\"ahler differentials as done in [L,\,Ch III,\,\S4]. However, here we follow ([M1,2]) to give a rather precise realization in 
terms of structures of two dimensional local fields.

\subsubsection{Residue maps for local fields}

\noindent
{\bf (A) Continuous differentials}
\vskip 0.20cm
\noindent
Let $(A,\m_A)$ be a local Noetherian ring and $N$ an $A$-module $N$. Denote by $N^{\sep}$ the maximal Hausdorff 
quotient of $N$ for the $\m_A$-adic topology, i.e., $N^\sep=N\big/\bigcap_{n=1}^\infty \m_A^nN$. In particular, if $A$ 
is an $R$-algebra for a certain ring $R$, then we have the differential module $\Omega_{A/R}$ and hence 
$\Omega_{A/R}^\sep$. Thus, if $F$ is a complete discrete valuation field and $K$ a subfield such that 
$\Frac(K\cap\O_F)=K$, then we have the space of the continuous differentials 
$$
\Omega_{F/K}^{\cts}:=  \Omega_{\O_F/\O_F\cap K}^\sep\,{\otimes}_{\O_F}F.
$$  
Consequently, if $F'/F$ is a finite, 
separable field extension, then $\Omega_{F'/K}^{\cts}=\Omega_{F/K}^{\cts}\otimes_FF'$ and hence there is a natural 
trace map $\Tr_{F'/F}: \Omega_{F'/K}^{\cts}\longrightarrow\ \Omega_{F/K}^{\cts}.$
\vskip 0.20cm
\noindent
{\bf (B) Equal characteristic zero}
\vskip 0.20cm
\noindent
Let $F$ be a two-dimensional local field of equal characteristic zero. Then $F$ contains a unique subfield $k_F$ of 
coefficients, up to isomorphism, such that $F\simeq k_F((t))$ for a suitable uniformizer $t$. In particular, 
$\Omega_{\O_F/\O_{k_F}}^\sep\simeq \O_F\cdot dt$ is a free $\O_F$-module of rank one. We define the residue 
map for $F$ by
$$
\res_F:\Omega_{F/k_F}^\cts\longrightarrow\ k_F,\qquad \omega=f\, dt\mapsto \coef_{t^{-1}}(f).
$$ 
By [M1], this is well defined, i.e., independent of the choice of $t$. Moreover, for a finite field extension $F'/F$, we 
have the following commutative diagram
$$
\begin{matrix}
\qquad\quad\Omega_{F'/k_F}^\cts&\xrightarrow{\ \res_{F'}\ }& k_{F'}\qquad\quad\\[0.40em]
\Tr_{F'/F}\downarrow&&\downarrow \Tr_{k_{F'}/k_F}\\[0.40em]
\qquad\quad\Omega_{F/k_F}^\cts&\xrightarrow{\ \res_{F}\ }& k_F.\qquad\quad
\end{matrix}
$$

\noindent
{\bf (C) Mixed characteristic}
\vskip 0.20cm
\noindent
Let $L$ be a two dimensional local field of mixed characteristics. Then the constant field $k_L$ of $L$ coincides with 
the algebraic closure of $\Q_p$ within $L$ for a certain prime number $p$, and $L$ itself is a finite field extension 
over $k_L\{\{t\}\}$ for a certain uniformizer $t$. Here, by definition,
$$
k_L\{\{t\}\}
:=\Big\{\sum_{i=-\infty}^\infty a_it^i\,:\,a_i\in k_L,\ \inf_i\big\{\nu_{k_L}(a_i)\big\}>-\infty\ (\forall i),\ \ \ 
a_i\to 0\ (i\to-\infty)\Big\}.
$$
Moreover, by [M1],  $\ \Omega_{\O_{k_L\{\{t\}\}/\O_{k_L}}}^\sep
=\O_{k_L\{\{t\}\}}\, dt\ \bigoplus\ \mathrm{Tors}\big(\Omega_{\O_{k_L\{\{t\}\}/\O_{k_L}}}^\sep\big).$ We\\[0.1em] 
define the residue map, first, for $k_L\{\{t\}\}$, by
$$
\res_{k_L\{\{t\}\}}: \Omega_{{k_L\{\{t\}\}/{k_L}}}^\cts\longrightarrow k_L,
\qquad \omega=f\, dt\mapsto -\coef_{t^{-1}}(f);
$$ 
then, for $L$, by the composition
$$
\res_L:
\Omega_{L/k_L}^\cts\xrightarrow{\,\Tr_{L/k_L\{\{t\}\}}\,}\,\Omega_{k_L\{\{t\}\}/k_L}^\cts
\xrightarrow{\,\res_{k_L\{\{t\}\}}\,} k_L.
$$
This is well defined by [M1]. Consequently, if $L'/L$ is a finite field extension, we have the commutative 
diagram
$$
\begin{matrix}
\qquad\quad\Omega_{L'/k_L}^\cts&\xrightarrow{\ \res_{L'}\ }& k_{L'}\qquad\quad\\[0.40em]
\Tr_{L'/L}\downarrow&&\downarrow \Tr_{k_{L'}/k_L}\\[0.40em]
\qquad\quad\Omega_{L/k_L}^\cts&\xrightarrow{\ \res_{L}\ }& k_L.\qquad\quad
\end{matrix}
$$

\subsubsection{Local residue maps}
As above, let $F$ be a number field and $\pi: X\to \Spec\,\O_F$ be an arithmetic surface with $X_F$ its generic fiber.

For each closed point $x\in X$, and a prime divisor $C$ on $X$ with $x\in C$,  by Theorem 1, the local ring $k(X)_{C,x}$ is a 
finite direct sum of two dimensional local fields, i.e., 
$$
k(X)_{C,x}=\bigoplus k(X)_{C_i,x},
$$ 
where $C_i$'s are normalized branches of the curve $C$ in a formal neighborhood $U$ of $x$. Set then 
$$
\res_{C,x}=\sum_i\res_{k(X)_{C_i,x}},
$$ 
which takes the values in $F_{\pi(x)}$, the local field of $F$ at the place $\pi(x)$. Recall also that, following [T], we 
have the canonical character 
$$
\lambda_{\pi(x)}:F_{\pi(x)}\xrightarrow{\Tr_{F_{\pi(x)}/\Q_p}}\Q_p\longrightarrow
\Q_p/\Z_p\longrightarrow\Q/\Z\hookrightarrow \R/\Z\simeq {\mathbb S}^1.
$$ 
Introduce accordingly 
$$
\Res_{C,x}:=\lambda_{\pi(x)}\circ \res_{C,x}.
$$

On the other hand, for each closed point $P\in X_F$,
$$
k(X_F)_P\,\widehat\bigotimes_F\,\Big(\prod_{\sigma\in S_\infty}F_\sigma\Big)
=\Big(\bigoplus\R((t))\Big)\bigoplus\Big(\bigoplus\mathbb C((t))\Big)
$$ 
is a finite direct sum of local fields $\R((t))$ and ${\mathbb C}((t)).$ Hence, similarly, for each $\sigma\in S_\infty$,
we have the associated residue maps $\res_{P,\sigma}$. Define 
$$
\Res_{P,\sigma}=\lambda_\sigma\circ\,\res_{P,\sigma}.
$$ 
Here, as in [T], to make all compatible, we set  $\displaystyle{\lambda_\sigma(x)=-{\Tr_{F_\sigma/\R}}}$, 
i.e., with a minus sign added.

\subsection{Global Residue Pairing}

The purpose here is to introduce a non-degenerate global residue pairing on the arithmetic adelic ring of an arithmetic 
surface.

\subsubsection{Global residue pairing}

Let $\pi:X\to\Spec\,\O_F$ be an arithmetic surface with $X_F$ its generic fiber. Then, by \S1.2.1, we have the 
associated arithmetic adelic ring 
$$
\A_X^{\ar}:=\A_{X,012}^{\ar}:=\A_X^{\fin}\,\bigoplus\, \A_X^\infty
$$ 
with $\A_X^{\fin}=\A_{X,012}$ the adelic ring for the 2-dimensional Noetherian scheme $X$ and 
$\A_X^\infty:=\A_{X_F}\widehat\otimes_{\mathbb Q}\mathbb R
:=\displaystyle{\lim_{\longrightarrow D_1}\lim_{\substack{\longleftarrow D_2\\ D_2\leq D_1}}
\Big(\big(\A_{X_F}(D_1)/\A_{X_F}(D_2)\big)\bigotimes_F\prod_{\sigma\in S_\infty}F_\sigma\Big).}$ By an abuse of 
notation, we will write elements of $\A_X^{\fin}$ as $(f_{C,x})_{C,x}$, or even $(f_{C,x})$,  and  elements of 
$\A_{X}^\infty$ as $(f_P)_{P}$, or even $(f_P)$.

Fix a rational differential $\omega=f(t)\, dt\not\equiv 0$ on $X$.
Then, we define a global pairing with respect to $\omega$ by
$$
\begin{matrix}
\langle\cdot,\cdot\rangle_\omega:\quad \A_X^\ar\times \A_X^\ar&\longrightarrow& \mathbb S^1\\
\big((f_{C,x},f_{P,\sigma}),(g_{C,x},g_{P,\sigma})\big)&\mapsto
&\sum_{C\subset X, x\in C:\,\pi(x)\in S_\fin}\Res_{C,x}(f_{C,x}g_{C,x}\omega)\\
&&\qquad\ \ +\sum_{P\in X_F}\sum_{\sigma\in S_\infty}\Res_{P,\sigma}(f_{P,\sigma}g_{P,\sigma}\omega).
\end{matrix}
$$

\begin{lem} Let $X$ be an arithmetic surface and $\omega$ is a non-zero rational differential on $X$. Then
the global pairing with respect to $\omega$ above is well defined.
\end{lem}
\noindent
{\it Proof.} Write 
$\sum_{C\subset X, x\in C:\,\pi(x)\in S_\fin}\Res_{C,x}(f_{C,x}g_{C,x}\omega)$ as a double summations
$\sum_{C\subset X}\sum_{ x\in C:\,\pi(x)\in S_\fin}\Res_{C,x}(f_{C,x}g_{C,x}\omega)$. 
Then, by Example 2, for all but finitely many curves $C$,  $\Res_{C,x}(f_{C,x}g_{C,x}\omega)=0$. So it suffices to 
show that for a fixed curve $C$, $\sum_{ x\in C:\,\pi(x)\in S_\fin}\Res_{C,x}(f_{C,x}g_{C,x}\omega)$ is finite. This is 
a direct consequence of the definition of $\A_{C,01}$ and  $\A_X^\fin$ in Example 2.

\subsubsection{Non-degeneracy}
\begin{prop} 
The residue pairing $\langle\cdot,\cdot\rangle_\omega$ on $\A_X^\ar$ is non-degenerate.
\end{prop}

\noindent
{\it Proof.} Let  $g\in \mathbb{A}_{X}^{\ar}$ be an adelic element such that, for all  $f\in \mathbb{A}_{X}^{\ar}$, 
$\langle f,g\rangle_\omega=0$. We show that $g=0$.

Rewrite the summation in the definition of $\langle\cdot,\cdot\rangle_\omega$ according to prime horizontal curves 
$E_P^\ar$ associated to closed points $P\in X_F$ and prime vertical curves $V\subset X$ appeared in the fibers of 
$\pi$. Namely, $\sum_{P\in X_F}\sum_{x\in E_P^\ar}+\sum_{V\in X}\sum_{x\in V}$. Then, note that, for a fixed adelic 
element $g=(g_{C,x},g_{P,\sigma})\in\A_X^\ar$,
$$
\langle (f_{C,x},f_{P,\sigma}),(g_{C,x},g_{P,\sigma})\rangle=0\qquad \forall (f_{C,x},f_{P,\sigma})\in\A_X^\ar.
$$ 
We have
$$
\sum_{P\in X_F}\sum_{x\in E_P^\ar}\Res_{E_P,x}(g_{E_P,x}f\omega)
+\sum_{V\in X}\sum_{x\in V}\Res_{E_P,x}(g_{E_P,x}f\omega)=0\qquad\forall f\in k(X).
$$ 

Now assume, otherwise, that $g\neq 0$. There exists either some vertical curve $C$ such that 
$0\neq g_{C, x}\in k(X)_{C, x}$, or, a certain algebraic point $P$ of the generic fiber $X_F$ such that  
$0\neq g_{P, \sigma}\in k(X_{F})_{P}\hat{\otimes}_{F}F_{\sigma}$. In case $g_{C, x}\neq 0$, by definition, 
$g_{C, x}=(g_{C_{i}, x})\in \bigoplus_i k(X)_{C_{i},x}=k(X)_{C, x}$, where $C_i$ runs over all branches of $C$, and 
$k(X)_{C_{i}, x}$ is a two-dimensional local field. Fix a branch $C_{i_0}$ such that $g_{C_{i_0},x}\neq 0$. By 
definition, ${\rm Res}_{C_{i}, x}:= \lambda_{\pi (x)}\circ {\rm res}_{C_{i}, x}$. So  we  can choose an element 
$h\in k(X)_{C_{i_0}, x}$ such that ${\rm Res}_{C_{i_0}, x}(h\omega)\neq 0$. For such $h$, we then take 
$f_{C_{i_0},x}\in k(X)_{C_{i_0}, x}$ such that $f_{C_{i_0}, x}g_{C_{i_0},x}=h$. Accordingly, if we construct an adelic 
element $f$ by taking all other components  $f_{C_{i}, x}$ to be zero but the $f_{C_{i_0}, x}$, then we have 
$\langle f, g\rangle_{\omega}={\rm Res}_{C_{i_0}, x}(f_{C_{i_0},x}g_{C_{i_0}, x}\omega)
={\rm Res}_{C_{i_0},x}(h\omega)\neq 0$. This contradicts to our original assumption that
$\langle f, g\rangle_{\omega}=0$. Hence, all the components of $g$ corresponding to vertical curves are zero.
Since we can use the same argument for the components of $g$ corresponding to algebraic points of $X_F$ to 
conclude that all the related components are zero as well, so $g=0$. This completes the proof.

\subsection{Adelic Subspaces}

\subsubsection{Level two subspaces}

Let  $\pi:X\to\Spec\,\O_F$ be an arithmetic surface. Our purpose here is to introduce certain level two intrinsic subspaces of 
$\A_X^{\ar}:=\A_{X,012}^{\ar}:=\A_X^{\fin}\oplus\A_X^\infty$. 

To start with, we analyze the structures of $\A_{X,01}^\fin$, one of the level two subspaces of 
$\A_X^{\fin}=\A_{X,012}^\fin$. By definition, see e.g.,  [P1], or better, 1.1.3(B), 
view as elements of $\mathbb A_X^{\mathrm{fin}}$ (via diagonal embedding), we may write
elements of $\mathbb A_{X,01}^{\mathrm{fin}}$ as $(f_{C,x})_{C,x}$. Then the partial components $(f_{C,x})_{x\in C}$ are independent of $x$, and hence can be written as $(f_C)_C$ (via diagonal embedding). 
Accordingly, we will simply write elements of $\A_{X,01}^\fin$ as $(f_C)_C$. On the other hand, with respect to $\pi$, curves 
on $X$ may be classified as being either vertical or horizontal. Therefore, we may and will write 
$(f_C)_C=(f_C)_{C:\,\mathrm{ver}}\times (f_C)_{C:\,\mathrm{hor}}.$ Accordingly, we set
$
\A_{X,01}^\fin=\A_{X,01}^{\fin,\mathrm{v}}\oplus \A_{X,01}^{\fin,\mathrm{h}},
$
where $\A_{X,01}^{\fin,\mathrm{v}}$, resp., $\A_{X,01}^{\fin,\mathrm{h}}$ denotes the collections of 
$(f_C)_{C:\,\mathrm{ver}}$, resp., $(f_C)_{C:\,\mathrm{hor}}.$ Furthermore, if $C$ is horizontal,  there exists an 
algebraic point $P$ of $X_F$ such that $C={\overline{\{P\}}}^X$, the Zariski closure of $P$ in $X$. For simplicity, 
write $C=E_P.$ Then $f_{E_P}\in F(X)_{E_P}$. But $F(X)_{E_P}=F(X_F)_P$. So it makes sense for us to talk about 
whether $f_{E_P}=f_P$ for a certain element $f_P\in F(X_F)_P$.

Fix a Weil divisor $D=D_v+D_h$ on $X$, where $D_v=\sum_Fn_VV$ with $V$ irreducible vertical curves 
and $D_h=\sum_Pn_PE_P$ with $E_P$ the horizontal curves. In 
particular, $D$ induces a divisor $D_F=\sum_Pn_PP$ on $X_F$. Following the uniformity condition in Definition 7(ii) (and that
for arithmetic curves recalled in \S1.2.4), we introduce  level two intrinsic  subspaces 
$\A_{X,01}^{\ar},\, \A_{X,02}^{\ar},\, \A_{X,12}^\ar(D)$  by 
$$
\begin{aligned} 
\A_{X,01}^{\ar}=&\big\{(f_{C,x})\times(f_P)\in\A_X^{\ar}\,\big|\,(f_{C,x})_{C,x}=(f_C)_{C,x}\in \A_{X,01}^{\fin},
\, f_{E_P}=f_P\ \forall P\in X_F\big\},\\[0.2em]
\A_{X,02}^{\ar}=&\A_{X,02}^\fin\bigoplus k(X_F)\widehat\otimes_{\mathbb Q}\R\quad\mathrm{and}\quad
\A_{X,12}^{\ar}(D):=\A_{X,12}^\fin(D)\bigoplus \big(\A_{X_F}(D_F)\widehat\otimes_{\mathbb Q}\R\big),
\end{aligned}$$
$$\begin{aligned}
\mathrm{where}\qquad\quad \A_{X,12}^\fin(D)
:=&\big\{(f_{C,x})\in \A_X^{\fin}\,\big|\, \ord_C(f_{C,x})+\ord_C(D)\geq 0\ \forall C\subset X\big\},\\[0.2em]
\A_{X_F}(D_F):=&\big\{(f_P)\in\A_{X_F} \,\big|\, \ord_P(f_P)+\ord_P(D_F)\geq 0\ \forall P\in X_F\big\},\\\end{aligned}$$
and 
$$
\A_{X_F}(D_F)\,\widehat\otimes_{\mathbb Q}\,\R
:=\lim_{\law D_F': D_F'\leq D_F}\Big(\A_{X_F}(D_F)/\A_{X_F}(D_F')\otimes_{\mathbb Q}\R\Big) .
$$
Here we have used the natural imbedding $k(X)= k(X_F)\hookrightarrow \A_{X_F}\hookrightarrow\A_{X}^{\ar}$.

Accordingly, we then also obtain three level one subspaces  
$$
\begin{aligned}
\A_{X,0}^{\ar}:=&\A_{X,01}^{\ar}\cap \A_{X,02}^{\ar},\qquad\qquad\mathrm{and}\\[0.30em]
\A_{X,1}^{\ar}(D):=&\A_{X,01}^{\ar}\cap \A_{X,12}^{\ar}(D),\qquad
\A_{X,2}^{\ar}(D):=\A_{X,02}^{\ar}\cap \A_{X,12}^{\ar}(D).
\end{aligned}
$$

\begin{lem} Let $X$ be an arithmetic surface, $D$ be a Weil divisor on $X$ with $D_F$ its induced divisor on $X_F$. 
We have
\vskip 0.150cm

\noindent
(i) $\A_{X,0}^{\ar}=k(X)$;
\vskip 0.150cm
\noindent
(ii) $\A_{X,1}^{\ar}(D)=\{(f_C)_{C,x}\times (f_P)\in \A_X^{\ar}: (f_C)_{C,x}\in \A_{X,1}(D),\ f_{E_P}=f_P\ \forall P\in X_F\}$;
\vskip 0.15cm
\noindent
(iii) $\A_{X,2}^{\ar}(D)=\{(f_x)_{C,x}\times (f)\in \A_X^{\ar}: \newline~\qquad \hskip 1.50cm 
(f_x)_{C,x}\in \A_{X,2}(D),\ f\in \A_{X_F}(D_F)\cap k(X_F)=H^0(X_F,D_F)\otimes_{\mathbb Q}\R\}$
\vskip 0.15cm
\noindent
(iv) Under the natural boundary map,  $(\A_\ar^*(X,D),d^*):$
$$
0\to  \A_{X,0}^\ar\oplus \A_{X,1}^\ar(D)\oplus \A_{X,2}^\ar(D)\to  
\A_{X,01}^\ar\oplus \A_{X,02}^\ar\oplus \A_{X,12}^\ar(D)\to\A_{X,012}^\ar\to 0
$$ 
forms a complex, the adelic complex for $D$.
\end{lem}

\noindent
{\it Proof.} 
The first three are direct consequences of the construction, while (iv) is standard from homotopy 
theory with the following the boundary maps: the first is the diagonal embedding, the 
second is given by $(x_0,x_1,x_2)\mapsto (x_0-x_1,x_1-x_2, x_2-x_0)$ and the final one is given by 
$(x_{01},x_{02},x_{12})$ $\mapsto x_{01}+x_{02}-x_{12}$. 

For example, 
to find  common elements $(f_{C,x}, f_P)$ of $\A_{X,01}^{\ar}$ and $\A_{X,02}^{\ar}$, 
we first concentrate on the infinite part.
Note that, by definition, the infinite part of  $\A_{X,02}^{\ar}$ is simply $k(X_F)\otimes \mathbb R$. 
Hence $(f_P)_P$ above, viewed as  elements of  infinite part of $\A_{X,01}^{\ar}$ should satisfy $f_P\in k(X_F)\otimes \mathbb R$. 
On the other hand, by the uniformity condition,
$f_P=f_{E_P}$\footnote{Here, as in our paper, to simplify notations,  this equality means that it holds in the related two dimensional local field after we embedd $k(X_F)$ as a subfield within.} with $f_{E_P}$ coming from finite part. This forces $f_P\in k(X_F)$ and hence $f_{E_P}\in k(X_F)$. That is to say $f_{{E_P} ,x}=f_P\in k(X_F)$. This, together with the fact that components of elements
in $\A_{X,01}^{\fin}$ (resp. $\A_{X,02}^{\fin}$) are independent of $x$ (resp. of $C$) then  completes the proof of (i). We leave others to the reader.

As a direct consequence of (ii) and (iii), we have the following

\begin{cor} Let $X$ be an arithmetic surface, $D$ be a Weil divisor on $X$. Then, we have the following
induced ind-pro structures on $\A_{X,01}^\ar$ and $\A_{X,02}^\ar$:
$$
\begin{aligned}
\A_{X,01}^\ar
=&\lim_{\raw D'}\lim_{\substack{\law D'\\D'\leq D}}\A_{X,1}^\ar(D)\big/\A_{X,1}^\ar(D'),\\
\A_{X,02}^\ar
=&\lim_{\raw D'}\lim_{\substack{\law D'\\D'\leq D}}\A_{X,2}^\ar(D)\big/\A_{X,2}^\ar(D').
\end{aligned}
$$
\end{cor}

Our definitions here are specializations of definitions in \S1.2.3 for $\FF=\O_X(D)$ over arithmetic surface $X$. 
We leave the details to the reader.

\subsubsection{Perpendicular subspaces}

For late use, we here establish a fundamental property for the level two arithmetic adelic subspaces introduced 
above. In fact, as we will see below, this property, in turn, characterizes these subspaces.
 
Fix a non-zero rational differential  $\omega$ on the arithmetic surface $X$. Then, by Proposition 12, we have a 
natural non-degenerate pairing $\langle\cdot,\cdot\rangle_\omega$ on $\A_X^{\ar}$. For a subspace $V$ of 
$\A_X^{\ar}$, set 
$$
V^\perp:=\big\{w\in\A_X^\ar\,\big|\ \langle w,v\rangle_\omega=0\ \ \forall v\in V\big\}
$$ 
be its perpendicular subspace of $V$ in $\A_X^\ar$ with respect to $\langle\cdot,\cdot\rangle_\omega$. Then, we 
have the following important

\begin{prop} 
Let $X$ be an arithmetic surface, $D$ be a Weil divisor and $\omega$ be a non-zero rational differential on $X$. 
Denote by $(\omega)$ the divisor on $X$ associated to $\omega$. Then we have

\noindent
(i)  $\qquad\Big(\A_{X,01}^{\ar}\Big)^\perp=\A_{X,01}^{\ar}$;

\noindent
(ii) $\quad\ \ \Big(\A_{X,02}^{\ar}\Big)^\perp=\A_{X,02}^{\ar}$;

\noindent
(iii)   $\quad\ \Big(\A_{X,12}^{\ar}(D)\Big)^\perp=\A_{X,12}^{\ar}((\omega)-D).$
\end{prop}
\noindent
{\it Proof.}  By an abuse of notation, we will write elements of $\A_X^{\fin}$ as $(f_{C,x})_{C,x}$ or even $(f_{C,x})$,  
and  elements of $\A_{X}^\infty$ as $(f_P)_{P}$ or even $(f_P)$.

We begin with a proof of $\A_{X,01}^{\ar}\subset (\A_{X,01}^{\ar})^\perp$. Let ${\bf f},\,{\bf g}\in \A_{X,01}^{\ar}$. By 
definition, ${\bf f}=(f_C)_{C,x}\times (f_P)_P$, ${\bf g}=(g_{C})_{C,x}\times (g_P)_P$, and, for all algebraic points $P$ 
of $X_F$, $f_{E_P}=f_P$ and $g_{E_P}=g_P$.  Consequently, 
$$
\begin{aligned}
\langle {\bf f},{\bf g}\rangle_\omega=&\sum_{x,C}\Res_{C,x}(f_Cg_C\omega)+\sum_{P}\Res_P(f_Pg_P\omega)\\
=&\sum_C\sum_{x:x\in C}\Res_{C,x}(f_Cg_C\omega)+\sum_P\Res_P(f_Pg_P\omega).
\end{aligned}
$$
If $C$ is a vertical curve on $X$, then, by the standard residue theorem for algebraic curves, see e.g. [P1], 
$\sum_{x:x\in C}\Res_{C,x}(f_Cg_C\omega)=0$. Hence
$$
\begin{aligned}
\langle {\bf f},{\bf g}\rangle_\omega 
=&\sum_{P\in X_F}\sum_{x\in E_P}\Res_P(f_{E_P}g_{E_P}\omega)+\sum_P\Res_P(f_Pg_P\omega)\\
=&\sum_{P\in X_F}\sum_{Q\in \overline E_P}\Res_Q(f_{P}g_{P}\omega).
\end{aligned}
$$ 
Here $\overline E_P$ denotes the Arakelov completion of $E_P$ associated to an algebraic point $P\in X_F$, and in the last
step, we have used our defining condition $f_{E_P}=f_P$ and $g_{E_P}=g_P$ for elements ${\bf f},\, {\bf g}$ in 
$\A_{X,01}^{\ar}$. Now, by  the residue theorem for  $\overline E_P$ ([M2, Thm 5.4]), 
$\sum_{Q\in \overline E_P}\Res_Q(f_{P}g_{P}\omega)=0$. Therefore, $\langle {\bf f},{\bf g}\rangle_\omega =0$, and 
$\A_{X,01}^{\ar}\subset (\A_{X,01}^{\ar})^\perp$. 
 
Next, we show that $\A_{X,01}^{\ar}\supset (\A_{X,01}^{\ar})^\perp$, based on the following ind-pro structure on 
$\A_{X,01}^{\ar}$ in Corollary 14:
$$
\A_{X,01}^\ar
=\lim_{\raw D'}\lim_{\substack{\law D'\\D'\leq D}}\A_{X,1}^\ar(D)\big/\A_{X,1}^\ar(D').
$$

Let $C$ be an irreducible curve $C$. Then, induced by the perfect pairing
$\langle\cdot,\cdot\rangle_\omega:\A_X^\ar\times\A_X^\ar\to\R/\Z$, for any divisor $D$ on $X$, by (iii), whose proof given 
below is independent of (i) and (ii), we obtain a pairing
$$
\A_{X,12}^{\ar}(D)\big/\A_{X,12}^{\ar}(D-C)\times \A_{X,12}^{\ar}((\omega)+C-D)\big/\A_{X,12}^{\ar}((\omega)-D)
\longrightarrow\R/\Z.\eqno(1)
$$ 
Moreover, directly from the definition, we have that $\A_{X,12}^{\ar}(D)\big/\A_{X,12}^{\ar}(D-C)\simeq \A_{C,01}^\ar$
for any divisor $D$. Hence, $\A_{X,12}^{\ar}((\omega)+C-D)\big/\A_{X,12}^{\ar}((\omega)-D)\simeq\A_{C,01}^\ar$ as 
well. So we can and will view $(1)$ as a pairing on $\A_{C,01}^\ar$.

If $C$ is vertical, then, there exists $\omega_C\in\Omega_{k(C)/\mathbb F_p}$ and an ${\bf a}=(a_v)\in\A^\ar_{C,01}$ 
such that (1) coincides with the pairing
$$
\langle\cdot,\cdot\rangle_{\omega_C,{\bf a}}:\A_{C,01}^\ar\times\A_{C,01}^\ar \longrightarrow\,\F_p;
\qquad ({\bf f},{\bf g})\mapsto{\sum}_v\Res_v(f_vg_va_v\omega_C).\eqno(2)
$$
Since $\A_{X,01}^{\ar}\subset (\A_{X,01}^{\ar})^\perp$, and, directly from the definition, we have that
$$
\A_{X,1}^{\ar}(D)\big/\A_{X,1}^{\ar}(D-C)\simeq \A_{X,1}^{\ar}((\omega)+C-D)\big/\A_{X,1}^{\ar}((\omega)-D)
\simeq\A_{C,0}^\ar=k(C),\eqno(3)
$$
we conclude that $\langle\cdot,\cdot\rangle_{\omega_C,{\bf a}}: \A_{C,01}^\ar\times \A^\ar_{C,01}\to\F_p$ annihilates 
$k(C)\times k(C)$. But $\langle\cdot,\cdot\rangle_{\omega_C,{\bf a}}$ can be identified with the canonical residue 
pairing $\langle\cdot,\cdot\rangle:\A_{C,01}^\ar\times \A_{C,01}^\ar\to\F_p$ associated to $C$. Consequently, 
$K(C)^\perp=K(C).\footnote{It is well-known that, see e.g., [Iw,\,\S4], if $\chi$ is a non-zero character on $\A_{C,01}^\ar$ 
such that $\chi(k(C))=\{0\}$, then the induced pairing 
$\langle\cdot,\cdot\rangle_\chi:\A_{C,01}^\ar\times \A_{C,01}^\ar\to\F_q;\, 
({\bf f},{\bf g})\mapsto \chi({\bf f}\cdot {\bf g})$ is perfect and $k(C)^\perp=k(C)$.}$ In particular, with respect to the 
pairing $\langle {\bf a}k(C),k(C)\rangle=0$. So ${\bf a}\in k(C)$. Hence, if necessary, with a possible modification on 
$\omega$, without loss of generality,  we may and will assume ${\bf a}=1$ and write 
$\langle\cdot,\cdot\rangle_{\omega_C,{\bf a}}$ simply as $\langle\cdot,\cdot\rangle_{\omega_C}$.
Therefore, by (3) and the fact that $k(C)^\perp=k(C)$, we have
$$
\Big(\A_{X,1}^{\ar}(D)\big/\A_{X,1}^{\ar}(D-C)\Big)^\perp\simeq 
\A_{X,1}^{\ar}((\omega)+C-D)\big/\A_{X,1}^{\ar}((\omega)-D).
$$
Moreover, with a verbatim change, the same discussion is valid for horizontal curves as well. Consequently, by 
applying this repeatedly, we have, for any irreducible curves $C_1, C_2$, the following commutating diagram with 
exact rows
$$
\begin{footnotesize}
\begin{matrix}
A(D-C_1)\big/A(D-C_1-C_2)&\hookrightarrow A(D)\big/A(D-C_1-C_2)\twoheadrightarrow &A(D)\big/A(D-C_1)\\[0.3em]
\|&\downarrow&\|\\[0.3em]
G\cap B(D-C_1)\big/G\cap B(D-C_1-C_2)&\hookrightarrow G\cap B(D)\big/G\cap  B(D-C_1-C_2)\twoheadrightarrow 
&G\cap B(D)\big/G\cap B(D-C_1)\\
\end{matrix}
\end{footnotesize}
$$
where, to save space, we set $A:=\A_{X,1}^\ar,\ B:=\A_{X,12}^\ar$ and $G:=\big(\A_{X,01}^\ar\big)^\perp$. 
Consequently, the vertical map in the middle is surjective. On the other hand, since 
$\A_{X,01}^\ar\subset \big(\A_{X,01}^\ar\big)^\perp$, this same map is also injective. Therefore, for any $D'\leq D$,
$$
\Big(\A_{X,1}^{\ar}(D)\big/\A_{X,1}^{\ar}(D')\Big)^\perp\simeq 
\A_{X,1}^{\ar}((\omega)-D')\big/\A_{X,1}^{\ar}((\omega)-D).
$$ 
with respect to our pairing
$$
\A_{X,1}^{\ar}(D)\big/\A_{X,1}^{\ar}(D')\times \A_{X,1}^{\ar}((\omega)-D')\big/\A_{X,1}^{\ar}((\omega)-D)
\longrightarrow\R/\Z.
$$
Consequently, we have
$$
\big(\A_{X,01}^\ar\big)^\perp=\A_{X,01}^\ar,
$$ 
since, by (iii) again, 
$$
\A_{X,01}^\ar
=\lim_{\raw D'}\lim_{\substack{\law D'\\D'\leq D}}\A_{X,1}^\ar(D)\big/\A_{X,1}^\ar(D')
=\lim_{\law D}\lim_{\substack{\raw D'\\D'\leq D}}\A_{X,1}^\ar((\omega)-D')\big/
\A_{X,1}^\ar((\omega)-D).
$$
This proves (i).

To prove (ii), we start with the inclusion  $\A_{X,02}^\ar\subset(\A_{X,02}^\ar)^\perp$. By definition, every element 
${\bf f}\in \A_{X,02}^\ar$ can be written as ${\bf f}=(f_x)_{C,x}\times(f_P)_P$ with $f_{C,x}=f_x$ and $(f_P)=(f)$ for 
some $f\in k(X_F)$ (since, by definition, the 02 type adeles are independent of one dimensional curves).
Thus, for ${\bf f},{\bf g}\in \A_{X,02}^\ar$, we have
$$
\begin{aligned}
\langle {\bf f},{\bf g}\rangle=&\sum_{C,x}\Res_{C,x}(f_{C,x}g_{C,x}\omega)+\sum_P\Res_P(f_Pg_P\omega)\\
=&\sum_x\sum_{C:C\ni x}\Res_{C,x}(f_xg_x\omega)+\sum_P\Res_P(fg\omega).
\end{aligned}
$$
Note that for a fixed $x$, $f_x$ and $g_x$ are fixed. Thus, by the residue theorem for the point $x$ ([M1, Thm 4.1]), 
we have $\sum_{C:C\ni x}\Res_{C,x}(f_xg_x\omega)=0$. So
$$
\langle {\bf f},{\bf g}\rangle=\sum_x0+\sum_P\Res_P(fg\omega)=\sum_P\Res_P(fg\omega).
$$
On the other hand, since $f,\,g\in k(X_F)$ and $\omega$ is a rational differential on $X_F$, the standard residue 
formula for the curve $X_F/F$ (see e.g., [S, \S II.7, Prop. 6]) implies that $\sum_P\Res_P(fg\omega)=0.$
Hence
$$
\langle {\bf f},{\bf g}\rangle=0,\qquad\forall \ {\bf f},{\bf g}\in \A_{X,02}^\ar.
$$
Therefore, $\A_{X,02}^\ar\subset(\A_{X,02}^\ar)^\perp$.

For the opposite direction $\A_{X,02}^\ar\supset(\A_{X,02}^\ar)^\perp$, similarly as in (i), we, in theory, can use the 
following ind-pro structure on $\A_{X,02}^\ar$:
$$
\A_{X,02}^\ar
=\lim_{\raw D'}\lim_{\substack{\law D'\\D'\leq D}}\A_{X,2}^\ar(D)\big/\A_{X,2}^\ar(D').
$$
However, due to the lack of details for horizontal differential theory in literature, we decide to first use this ind-pro 
structure to merely prove the part that if ${\bf f}=(f_{C,x})_{C,x}\times (f_P)_P \in (\A_{X,02}^\ar)^\perp$,  for vertical 
$C$'s, $f_{C,x}=f_x$; and then to take a more classical approach for the rest.

Choose ${\bf f}=(f_{C,x})_{C,x}\times (f_P)_P\in (\A_{X,02}^\ar)^\perp$. Then by our assumption, for any element  
${\bf g}=(g_x)_{C,x}\times(g)_P\in \A_{X,02}^\ar,$ we have
$$
0=\langle {\bf f},{\bf g}\rangle=\sum_{x}\sum_{C:C\ni x}\Res_{C,x}(f_{C,x}g_{x}\omega)+\sum_P\Res_P(f_Pg\omega).
\eqno(4)
$$
Now note that the element $(g_x)_{C,x}\in\A_{X,02}^\ar$ and the element $g\in k(X_F)$ can be changed totally 
independently, we conclude that both of the summations, $\sum_{x}\sum_{C:C\ni x}\Res_{C,x}(f_{C,x}g_{x}\omega)$ 
and $\sum_P\Res_P(f_Pg\omega)$, are constants independent of ${\bf g}$. This then implies that both of them are 0. 
Indeed, since $\sum_P\Res_P(f_Pg\omega)$ is a constant independent of $g$, by choosing $g$ to be constant 
function, we conclude that $\ \sum_P\Res_P(f_Pg\omega)\equiv0,\ \, \forall g\in k(X_K).$ Consequently, by (4), we 
have, from (2),
$$
\sum_{x}\sum_{C:C\ni x}\Res_{C,x}(f_{C,x}g_{x}\omega)\equiv 0,\ \ \forall (g_x)\in A_{X,02}^\ar.\eqno(5)
$$

To end the proof, we need to show that  $f_{C,x}$'s are independent of $C$ and $f_P$ are independent of $P$. 
First, we treat the case when $C$ is vertical. As said above, we will use the associated ind-pro structures. So, 
assume for now that $C$ is vertical. Then, for any divisor $D$ on $X$, we have 
$$
\begin{aligned}
\A_{X,2}^\ar(D)\big/\A_{X,2}^\ar(D-C)\simeq&\A_{C,1}(D|_C),\\[0.30em]
\A_{X,2}^\ar((\omega)+C-D)\big/\A_{X,2}^\ar((\omega)+C-(D-C))\simeq&\A_{C,1}((\omega_C')-D|_C),
\end{aligned}
$$
for a certain $\omega_C'\in\Omega_{k(C)/\F_p}$ satisfying $(\omega_C')=\big((\omega)+C\big)|_C$ by the adjunction formula.
We claim that $(\omega_C')=(\omega_C)$. Indeed, since $\A_{C,1}(D|_C)^\perp=\A_{C,1}((\omega_C)-D|_C)$ and 
$\A_{X,02}^\ar\subset\big(\A_{X,02}^\ar\big)^\perp$, $\A_{C,1}((\omega_C')-D|_C)\subset \A_{C,1}((\omega_C)-D|_C)$. This 
implies that $(\omega_C)\geq (\omega_C')$ and hence $(\omega_C)= (\omega_C')$, because there is no $f\in k(C)$ such that 
$(f)>0$. Thus, with respect to the canonical residue pairing on $C$, we have 
$\A_{C,1}(D|_C)^\perp=\A_{C,1}((\omega_C')-D|_C)$. Consequently, as in (i), from 
$$
\A_{X,02}^\ar
=\lim_{\raw D'}\lim_{\substack{\law D'\\D'\leq D}}\A_{X,2}^\ar(D)\big/\A_{X,2}^\ar(D')
=\lim_{\law D}\lim_{\substack{\raw D'\\D'\leq D}}
\A_{X,2}^\ar((\omega)-D')\big/\A_{X,2}^\ar((\omega)-D),
$$
we conclude that $f_{C,x}=f_x\in k(X)_x$ and hence independent of $C$.

To prove the rest, we take a classical approach with a  use of Chinese reminder theorem, using an idea in the proof of 
Proposition 1 of [P1]. To be more precise, fix $x_0\in X$ and  a prime divisor $H$ on $X$ satisfying  that $x_0\in H$ and  that 
$V:=\mathrm{Spec}\,\O_{X,x_0}-H$ is affine. We claim that for any family of prime divisors $D_j$ on $V$ such that 
$D_i\cap D_j=\emptyset$ if $i\not =j$, and any rational functions $f_0,f_1,\dots, f_n$ on $V$, and any fixed divisor $D$ 
supported on $D_i$'s, there exists a rational function $g$ such that
$$
\begin{cases}
\ord_{D_i}(f_i-g)\geq\ord_{D_i}(D),& i=0,1,\dots, n,\\
\ord_{D_i}(g)\geq\ord_{D_i}(D),& i\not\in\{0,1,\dots,n\}
\end{cases}
$$
Indeed, by clearing the common denominators for $f_i$'s,  (with a modification of $f_i$'s if necessary,) we may 
assume that $f_i$'s are all regular. Then by applying the Chinese reminder theorem to the fractional ideals
$\frak p_i^{\ord_{D_i}(D)}, i=0,1,\dots, n$, and $\cap_{i\not\in\{0,1,\dots,n\}}\frak p_i^{\ord_{D_i}(D)}$, where 
$\frak p_i$ are the prime ideas associated to the prime divisors $D_i$, we see the existence of such a $g$. 

Associated to ${\bf f}\in \Big(\A_{X,02}^\ar\big)^\perp$, form a new adele ${\bf f}'\in\A_{X,012}^\ar$ by setting 
$
f_{C,x}'
=f_{C,x}-f_{H,x}$ where $H$ is a fixed vertical curve. Since, as proved above, for vertical $H$, $f_{H,x}=f_x\in k(X)_x$,
this is well defined.
Then 
$$
\begin{aligned}
&\sum_{C:C\ni x_0}\Res_{C,x_0}(f_{C,x_0}'g_{x_0}\omega)\\
=&\sum_{C:C\ni x_0}\Res_{C,x_0}(f_{C,x_0}g_{x_0}\omega)
-\sum_{C:C\ni x_0}\Res_{C,x_0}(f_{H,x_0}g_{x_0}\omega).
\end{aligned}
$$
The first sum is zero, since we have (5) by our choice of ${\bf f}$. The second sum, being taken over all prime curves 
passing through $x_0$, is zero as well, since we can apply the residue theorem for the point $x_0$ as above (with 
$f_{H,x_0}$ being independent of $C$). That is to say, 
$$
\sum_{C:C\ni x_0}\Res_{C,x_0}(f_{C,x_0}'g_{x_0}\omega)=0.
$$

Now, applying the above existence to obtain a $g$ satisfying that for any fixed rational function $f_0$ and any fixed 
curve $C_0\ni x_0$, we have
$$
\begin{cases}
\ord_{C}(f_{C,x_0}')+\ord_{C}(f_0-g)+\ord_{C}(\omega)\geq 0,&C=C_0\\
\quad \ord_{C}(f_{C,x_0}')+\ord_{C}(g)+\ord_{C}(\omega)\geq 0,& C\not=C_0, H.
\end{cases}
$$
Consequently, by the definition of the residue map, with $f_{H,x_0}'\equiv 0$ in mind, we get, for any $f_0\in k(X)$ and 
the corresponding $g$ just chosen, 
$$
\begin{aligned}
0=&\sum_{C:C\ni x_0}\Res_{C,x_0}(f_{C,x_0}'g\omega)
=\sum_{C:C\ni x_0,\{C:C\not=C_0,H\}\cup \{C:C=C_0,H\}}\Res_{C,x_0}(f_{C,x_0}'g\omega)\\
=&\sum_{C:C\ni x_0,C=C_0,H}\Res_{C,x_0}(f_{C,x_0}'g\omega)=\Res_{C_0,x_0}(f_{C_0,x_0}'g\omega)
=\Res_{C_0,x_0}(f_{C_0,x_0}'f_0\omega).
\end{aligned}
$$ 
Since the last quantity is always zero for all $f_0$, this then implies that $f_{C,x_0}'=0$, namely, 
$f_{C_0,x_0}=f_{H,x_0}$. 

To end the proof of (ii), we still need to show that $f_P=f_{P_0}$ for a fixed $P_0\in X_K$ and all $P\in X_K$. But this 
is amount to a use of a similar argument just said again, based on Chinese reminder theorem. See e.g., [Iw, \S 4]. We leave 
the details for the reader. Thus, if ${\bf f}=(f_{C,x})_{C,x}\times (f_P)_P\in (\A_{X,02}^\ar)^\perp$, then 
$f_{C,x}=f_{C_0,x}$ and $f_P=f_{P_0}$ for fixed $C_0$ and $P_0$. Therefore, ${\bf f}\in \A_{X,02}^\ar.$
That is to say, $(\A_{X,02}^\ar)^\perp=\A_{X,02}^\ar.$ This proves (ii).
\vskip 0.20cm
Finally, we prove (iii). The inclusion $\displaystyle{\A_{X, 12}^{\ar}((\omega)-D)\subset (\A_{X, 12}^{\ar}(D))^{\bot}}$ 
is easy.\\[0.20em]
Indeed, for $\displaystyle{{\bf f}=(f_{C, x})\times (f_{P})\in \A_{X,12}^{\ar}((\omega)-D)}$, 
$\displaystyle{{\bf g}=(g_{C, x})\times (g_{P})\in \A_{X,12}^{\ar}(D),}$\\[0.20em] we have $\ \langle {\bf f}, {\bf g}\rangle
={\sum}_{C, x}{\rm Res}_{C, x}(f_{C, x}g_{C, x}\omega)+{\sum}_{P}{\rm Res}_{P}(f_{P}g_{P}\omega)=0,\ $
since every\\[0.20em] 
term in each of these two summations is zero by definition.

To prove the other direction $\A_{X, 12}^{\ar}((\omega)-D)\supset (\A_{X, 12}^{\ar}(D))^{\bot}$, we make the following 
preparations. Set $\widetilde{\pi}:X\rightarrow {\rm Spec}\,{\cal O}_{F}\rightarrow {\rm Spec}\,\mathbb{Z}$. Then,  by 
[M1, Thm 5.7], the dualizing sheaf $\omega_{\widetilde{\pi}}$ of $\widetilde{\pi}$ is given by, for an open subset 
$U\subseteq X$, 
$$
\omega_{\widetilde{\pi}}(U)=\big\{\omega\in \Omega_{k(X)/\mathbb{Q}}\mid {\rm Res}_{C, x}(f\omega)=0\quad 
\forall  x\in C(\subset U),\ \forall  f\in {\cal O}_{X, C}\big\}.
$$
By a similar argument as in [M1] (used to prove the above result), we have, for a fixed curve $C_{0}$,
$$
\omega_{\widetilde{\pi}, C_{0}}=\big\{\omega\in \Omega_{k(X)/\mathbb{Q}}\mid {\rm Res}_{C_{0}, x}(f\omega)=0\quad
\forall \,x\in C_{0},\ f\in {\cal O}_{X, C_{0}}\big\}.
$$
This is nothing but the collection of differentials $\omega$ satisfying ${\rm ord}_{C_{0}}((\omega))\geq 0$. Moreover, 
we have, for a fixed pair $x_{0}\in C_{0}$,
$$
\omega_{\widetilde{\pi}, C_{0}}\otimes_{{\cal O}_{C_{0}}}{\cal O}_{C_{0}, x_{0}}
=\big\{\omega\in \Omega_{k(X)_{C_{0}, x_{0}}/\mathbb{Q}_{\widetilde{\pi}(x_{0})}}^{{\rm cts}}\mid 
{\rm Res}_{C_{0}, x_{0}}(f\omega)=0\quad{\forall}\,f\in {\cal O}_{C_{0}, x_{0}}\big\}.
$$ 
This is simply the collection of differentials $\omega$ satisfying ${\rm ord}_{C_{0}}((\omega))\geq 0$. From these, we 
conclude that the following conditions are equivalent for a non-zero differential 
$\omega\in \Omega_{k(X)/\mathbb{Q}}$ and a fixed pair $x_{0}\in C_{0}$;\\[0.2em]
(1)$_C$\ \ $\forall f\in {\cal O}_{X, C_{0}}\ \ {\rm Res_{C_{0}, x_{0}}}(f \omega)=0$.\\[0.2em]
(2)$_C$\ \ ${\rm ord}_{C_{0}}((\omega))\geq 0.$\\
Furthermore, a similar argument for points instead of curves also gives us that the following conditions are equivalent 
for a fixed point $P$;\\[0.2em]
(1)$_P$\ \ $\forall f\in {\cal O}_{X_{F}, P_{0}}\ \ {\rm Res}_{P_{0}}(f \omega)=0$.\\[0.2em]
(2)$_P$\ \ ${\rm ord}_{P_{0}}((\omega))\geq 0.$

Now we are ready to continue our proof. Let ${\bf f}=(f_{C, x})\times (f_{P})\in \A_{X,12}^{\ar}(D)^{\bot}$ and 
${\bf g}=(g_{C, x})\times (g_{P})\in \A_{X,12}^{\ar}(D)$. Then
$$
0=\langle {\bf f},{\bf g}\rangle
={\sum}_{C, x}{\rm Res}_{C, x}(f_{C, x}g_{C, x}\omega)+{\sum}_{P}{\rm Res}_{P}(f_{P}g_{P}\omega).
$$
Thus, by the independence of components of adeles, we see that, for fixed $(C, x)$ and $P$,
${\rm Res}_{C, x}(f_{C, x}g_{C, x}\omega)=0$ and ${\rm Res}_{P}(f_{P}g_{P}\omega)=0$. This, together with the 
equivalence between (1)$_C$ and (2)$_C$ and the equivalence between (1)$_P$ and (2)$_P$, we conclude that 
${\bf f}=(f_{C, x})\times (f_{P})\in \A_{X, 12}^{\ar}((\omega)-D)$, and hence complete the proof.

\subsection{Arithmetic Cohomology Groups}

\subsubsection{Definitions}

Let $X$ be an arithmetic surface and $D$ be a Weil divisor on $X$. By Lemma 13, particularly, with  (iv), we obtain 
an arithmetic adelic complex.  Accordingly, we define the {\it arithmetic cohomology groups $H_\ar^i(X,D)$} by 
$H^i(\A_\ar^*(X,D),d^*)$, the {\it i-th} cohomology groups of the complex $(\A_\ar^*(X,D),d^*)$, $i=0,1,2$. 

\begin{prop}  (i) The arithmetic cohomology groups $H_\ar^i(X,D)$ of  $D$ on $X$, $i=0,1,2$, are given by
$$\begin{aligned}H^0_{\ar}(X,D)=&\A_{X, 01}^{\ar}\cap \A_{X, 02}^{\ar}\cap \A_{X, 12}^{\ar}(D);\\
H^1_{\ar}(X,D) =&\\
=\Big(\big(\A_{X, 01}^{\ar}+&\A_{X, 02}^{\ar}\big)\cap \A_{X, 12}^{\ar}(D)\Big)\Big/
\Big(\A_{X, 01}^{\ar}\cap \A_{X, 12}^{\ar}(D)+\A_{X, 02}^{\ar}\cap \A_{X, 12}^{\ar}(D)\Big);\\
H^2_{\ar}(X,D)=&\A_{X, 012}^{\ar}\Big/\Big(\A_{X, 01}^{\ar}+ \A_{X, 02}^{\ar}+\A_{X, 12}^{\ar}(D)\Big).\end{aligned}$$
(ii) There exist natural isomorphisms
$$\begin{aligned}
&H^1_{\ar}(X,D)\\
\simeq&\Big(\A_{X, 01}^{\ar}\cap \big(\A_{X, 02}^{\ar}+ \A_{X, 12}^{\ar}(D)\big)\Big)\Big/
\Big(\A_{X, 01}^{\ar}\cap \A_{X, 02}^{\ar}+\A_{X, 01}^{\ar}\cap \A_{X, 12}^{\ar}(D)\Big)\\
\simeq&\Big(\A_{X, 02}^{\ar}\cap \big(\A_{X, 01}^{\ar}+\A_{X, 12}^{\ar}(D)\big)\Big)\Big/\Big(\A_{X, 01}^{\ar}\cap \A_{X, 02}^{\ar}
+\A_{X, 02}^{\ar}\cap \A_{X, 12}^{\ar}(D)\Big)
\end{aligned}$$
\end{prop}

\noindent
{\it Proof.} This is an arithmetic analogue of [P1]. Indeed, (i) comes directly from the definition, and via the first and second isomorphism
theorems in elementary group theory, (ii) is obtained using
a direct group theoretic calculation for
our arithmetic  adelic complex.

As mentioned above, the cohomology groups for the divisor $D$ defined here coincide with the adelic global cohomology 
defined in \S 1.2.3 associated to the invertible sheaf  $\O_X(D)$. That is,  $H^i_{\ar}(X,D)=H^i_{\ar}(X,\mathcal O_X(D))$. In 
fact, our general construction in \S 1.2.3 for quasi-coherent sheaves over higher dimensional arithmetic varieties is obtained 
modeling our constructions for arithmetic curves and arithmetic surfaces.

\subsubsection{Inductive long exact sequences}

\noindent
{\bf (V) Vertical Curves}

Just like the classical cohomology theory,  arithmetic cohomology groups also admit an inductive structure related to 
vertical geometric curves  on arithmetic surfaces. More precisely, we have the following

\begin{prop} 
Let $C$ be an irreducible  vertical curve of $X$, then, for any  $D$,  we have the long exact sequence of cohomology 
groups
$$
\begin{aligned}
0\to &H^0_{\ar}(X,C)\to H^0_{\ar}(X,D+C)\to H^0_{\ar}(C,(D+C)|_C)\\
&\to H^1_{\ar}(X,D)\to H^1_{\ar}(X,D+C)\to H^1_{\ar}(C,(D+C)|_C)\\
&\to H^2_{\ar}(X,D)\to H^2_{\ar}(X,D+C)\to 0.
\end{aligned}
$$
Here $H^i_{\ar}(C,(D+C)|_C),\, i=0,1$ are the usual cohomology groups  for vertical geometric curves.
\end{prop}
\noindent
{\it Proof.}  For $C$ be an irreducible vertical curve in $X$. Then from definition, one calculates that

\noindent
(a) $\dis{\A_{X,12}^{\ar}(D+C)/A_{X,12}^{\ar}(D)
=\A_{X,12}^{\fin}(D+C)/A_{X,12}^{\fin}(D)\oplus\{0\}=\A_{C,01}\oplus \{0\}.}$

\noindent
(b) $\A_{X,1}^{\ar}(D+C)/\A_{X,1}^{\ar}(D)=k(C)$ as there are neither changes along horizontal curves nor along $X_F$.

\noindent
(c) $\A_{X,2}^{\ar}(D+C)/\A_{X,2}^{\ar}(D)$
$$
=\A_{X,2}^{\fin}(D_\fin+C)/\A_{X,2}^{\fin}(D_\fin)\oplus\{0\}=\A_C\Big((D_\fin+C)|_{C}\Big)\oplus\{0\}.
$$

Consequently, for the morphism
$$
\begin{matrix}
\A_{X,1}^{\ar}(D+C)/\A_{X,1}^{\ar}(D)\oplus \A_{X,2}^{\ar}(D+C)/\A_{X,2}^{\ar}(D)&\buildrel\phi\over \longrightarrow&
\A_{X,12}^{\ar}(D+C)/A_{X,12}^{\ar}(D),\\  (x,y)&\mapsto &x-y,
\end{matrix}
$$ 
we conclude that

\noindent
(d) $\mathrm{Ker}\,\phi$ is given by
$$
\Big(\A_{X,1}^{\ar}(D+C)/\A_{X,1}^{\ar}(D)\Big)\cap \Big(\A_{X,2}^{\ar}(D+C)/\A_{X,2}^{\ar}(D)\Big)
=H^0_{\ar}(C,(D+C)|_C),
$$

\noindent
(e) $\mathrm{Coker}\,\phi$ is given by
$$
\begin{aligned}
\Big(\A_{X,12}^{\ar}&(D+C)/A_{X,12}^{\ar}(D)\Big)
\Big/ \Big(\A_{X,12}^{\ar}(D+C)/A_{X,12}^{\ar}(D)+\A_{X,12}^{\ar}(D+C)/A_{X,12}^{\ar}(D)\Big)\\
=&H^1_{\ar}(C,(D+C)|_C).
\end{aligned}
$$
Therefore, by definition, we have the long exact sequence
$$
\begin{aligned}
0\to &H^0_{\ar}(X,C)\to H^0_{\ar}(X,D+C)\to H^0_{\ar}(C,(D+C)|_C)\\
&\to H^1_{\ar}(X,D)\to H^1_{\ar}(X,D+C)\to H^1_{\ar}(C,(D+C)|_C)\\
&\to H^2_{\ar}(X,D)\to H^2_{\ar}(X,D+C)\to 0.
\end{aligned}
$$

This then completes the proof.

\vskip 0.30cm
\noindent
{\bf (H) Horizontal Curves}

For horizontal curve $E_P$ corresponding to an algebraic point $P$ of $X_F$,  we have

\noindent
(a) $\A_{X,12}^{\ar}(D+E_P)/A_{X,12}^{\ar}(D)$
$$
\begin{aligned}
=&\A_{X,12}^{\fin}(D+E_P)/A_{X,12}^{\fin}(D)\oplus\A_{X_F}(D_F+P)/A_{X_F}(D_F)\otimes_{\mathbb{Q}}\R\\
=&\A_{E_P,01}\oplus k(E_P)\otimes_{\mathbb{Q}}\R=\A_{E_P}^{\ar};
\end{aligned}
$$

\noindent
(b) $\A_{X,1}^{\ar}(D+E_P)/\A_{X,1}^{\ar}(D)=k(E_P)$,  diagonally embedded in
$$
\A_{X,1}^{\fin}(D_\fin+E_P)/\A_{X,1}^{\fin}(D_\fin)\oplus\A_{X_F}(D_F+P)/\A_{X_F}(D_F)\otimes_{\mathbb{Q}}\R;
$$

\noindent
(c) $\A_{X,2}^{\ar}(D+E_P)/\A_{X,2}^{\ar}(D)$
$$
=\A_{E_P}\Big((D_\fin+E_P)|_{E_P}\Big)\oplus H^0(X_F,D_F+E_P)\otimes \R/H^0(X_F,D_F)\otimes\R.
$$

Similarly, we have the corresponding morphism
$$
\begin{matrix}
\A_{X,1}^{\ar}(D+E_P)/\A_{X,1}^{\ar}(D)\oplus \A_{X,2}^{\ar}(D+E_P)/\A_{X,2}^{\ar}(D)
&\buildrel\varphi\over \longrightarrow&\A_{X,12}^{\ar}(D+E_P)/A_{X,12}^{\ar}(D)\\ 
(x,y)&\mapsto &x-y,
\end{matrix}
$$ 
and hence obtain the following proposition in parallel.

\begin{prop} 
Let $E_P$ be the horizontal curve $E_P$ corresponding to an algebraic point $P$ of $X_F$, then, for any  $D$, We have the long exact sequence
$$
\begin{aligned}
0\to &H^0_{\ar}(X,E_P)\to H^0_{\ar}(X,D+E_P)\to \mathrm{Ker}\,\varphi\\
&\to H^1_{\ar}(X,D)\to H^1_{\ar}(X,D+E_P)\to \mathrm{Coker}\,\varphi\\
&\to H^2_{\ar}(X,D)\to H^2_{\ar}(X,D+E_P)\to 0.
\end{aligned}
$$
\end{prop}
However, unlike for vertical curves, we do not have the group isomorphisms between
$\mathrm{Ker}\,\varphi$, resp. $\mathrm{Coker}\,\varphi$, and $H^0_{\ar}(E_P,(D+E_P)|_{E_P})$, resp.  
$H^1_{\ar}(E_P,(D+E_P)|_{E_P})$. This is in fact not surprising: different from vertical curves, for the arithmetic 
cohomology, there is no simple additive law with respect to horizontal curves when count these arithmetic groups: In 
Arakelov theory, we only have
$$
\chi_{\ar}(X,D+E_P)=\chi_{\ar}(X,D)+\chi_{\ar}(E_P, (D+E_P)|_{E_P})-\frac{1}{2}d_\lambda(E)\eqno(6)
$$
with  discrepancy $-\frac{1}{2}d_\lambda(E)$ resulting from Green's functions. (See e.g., [L, p.114].)

On the other hand, recall that, on generic fiber $X_F$, we have  the long exact sequence of cohomology groups
$$
\begin{aligned}
0\to& H^0(X_F,D_F)\to H^0(X_F, D_F+P)\to \O_{X_F}(D+P)|_P\to Q\to 0\\
&(\text{with }Q\text{ defined by}\qquad 0\to Q\to H^1(X_F,D_F)\to H^1(X_F, D_F+P)\to 0),
\end{aligned}\eqno(7)
$$
and that, in Arakelov theory, see e.g., [L, VI, particularly, p.140], what really used is the much rough version 
$\lambda(D_F+P)\simeq \lambda(D_F)\otimes \O_{X_F}(D+P)|_P$ where $\lambda$ denotes the 
Grothendieck-Mumford determinant. One checks that the exact sequence (6) does appear in our calculation above. 
Indeed, for curve $X_F/F$, 
$$
H^0(X_F,D_F)=k(X_F)\cap \A_{X_F}(D_F)\quad\&\quad H^1(X_F,D_F)=\A_{X_F}/k(X_F)+ \A_{X_F}(D_F).
$$
Consequently, (7) is equivalent to the exact sequence
$$
\begin{aligned}
0\to  H^0(X_F,D_F+P)/&H^0(X_F,D_F)\to \A_{X_F}(D_F+P)/A_{X_F}(D_F)\to Q\to 0\\
&\hskip 1.0cm 0\to Q\to  H^1(X_F,D_F)\to H^1(X_F,D_F+P)\to 0.
\end{aligned}
$$ 
(Note that $\A_{X_F}(D_F+P)/A_{X_F}(D_F)$ is supported only on $P$.) Clearly, all this can be read from the 
calculations in (a,b,c) and the morphism $\varphi$ above. So our construction offers a much more refined structure 
topologically.

\subsubsection{Duality of cohomology groups}

Let $\pi:X\to \Spec\,\O_F$ be an arithmetic surface defined over  the ring of integers of a number field $F$. Then we 
have the adelic space $\A_X^{\ar}$, its level two subspaces $\A_{01}^{\ar}, \,\A_{02}^{\ar}$ and $\A_{12}^{\ar}(D)$,  and 
hence the cohomology groups $H_\ar^i(X,\O_X(D))$ associated to a Weil divisor $D$ on $X$.  Moreover, there is a 
natural ind-pro structure on $\A_X^\ar$
$$
\A_X^\ar=\lim_{\raw D}\lim_{\law E: E\leq D}\A_{X,12}(D)\big/\A_{X,12}(E).
$$ 
Note that $\A_{X,12}(D)\big/\A_{X,12}(E)$'s are locally compact topological spaces. Consequently, as explained in \S 3, the 
next section, induced from the projective limit, we get a natural final topology on 
$\dis{\A_{X,12}(D)=\lim_{\law E: E\leq D}\A_{X,12}(D)\big/\A_{X,12}(E)}$; similarly, induced from the inductive limit, we 
get a natural initial topology on $\dis{\A_X^\ar=\lim_{\raw D}\A_{X,12}(D)}$. Moreover, by Theorem II, we know that
all three level two subspaces  $\A_{01}^{\ar}, \,\A_{02}^{\ar}$ and $\A_{12}^{\ar}(D)$ are closed in  $\A_X^\ar$.
Consequently, we obtain natural topological structures on arithmetic cohomology groups $H_\ar^i(X,\O_X(D))$ induced 
from the canonical topology of $\A_X^\ar$. Our main theorem here is the following

\begin{thm}  Let $X$ be an arithmetic surface and $D$ be a Weil divisor on $X$. Then as topological groups,
we have natural isomorphisms
$$
\widehat {H^i_{\ar}(X,D)}\simeq H^{2-i}_{\ar}((\omega)-D)\qquad i=0,1,2.
$$
\end{thm}
Recall that, for a topology space $T$,  its topological dual is defined by
$\widehat{\,T\,}:=\{\phi:T\to \mathbb S^1\ \mathrm{continuous}\}$ together with an compact-open topology. (See e.g., \S 3.1.1 
for details.)

This theorem is proved at the end of this paper, after we expose some basic structural results for the ind-pro topology on 
$\A_X^\ar$ in the next  section.

\eject
\section{Ind-Pro Topology in Dimension Two}

In this section, we establish some basic properties for  ind-pro topologies on various adelic spaces associated to 
arithmetic surfaces. This may be viewed as a natural generalization of a well-known topological theory for one dimensional 
adeles (see e.g. [Iw], [T]). Our main result here is the following
\vskip 0.20cm
\noindent
{\bf Theorem II.}
{\it Let $X$ be an arithmetic surface. Then with respect to the canonical ind-pro topology on $\A_X^\ar$, we have

\noindent
(1) Level two subspaces $\A_{X,01}^\ar, \A_{X,02}^\ar$ and $\A_{X,12}^\ar(D)$ are closed in $\A_X^\ar$;

\noindent
(2) $\A_X^\ar$ is a Hausdorff, complete, and compact oriented topological group;

\noindent
(3) $\A_X^\ar$ is self-dual. That is, as topological groups,}
$$
\widehat{\ \A_X^\ar\ }\ \simeq\ \A_X^\ar.
$$

\subsection{Ind-pro topologies on adelic spaces}

\subsubsection{Ind-pro topological spaces and their duals}

To begin with, let us recall some basic topological constructions for inductive limits and projective limits of topological 
spaces. 

\noindent
(1) Let $\{G_m\}_m$ be an inductive system of topological spaces, $\displaystyle{G:=\lim_{\ \longrightarrow m}G_m}$ 
with structure maps $\iota_m:G_m\to G$. Then, the inductive topology on $G$ is defined by assigning 
subsets $U$ of $G$ to be open, if $\iota_m^{-1}(U)$ is open in $G_m$ for each $m$. Inductive topology is also 
called the final topology since it is the finest topology on $G$ such that $\iota_m:G_m\to G$ are continuous.

\noindent
(2) Let $\{G_n\}_n$ be a projective system of topological spaces, $\displaystyle{G:=\lim_{\longleftarrow n}G_n}$ with 
structure maps $\pi_n:G\to G_n$. Then,  the projective topology on $G$ is defined as  the one generated by open 
subsets $\pi_n^{-1}(U_n)$, where $U_n$ are open subsets of $G_n$. Projective topology is also called the initial 
topology since it is the coarsest topology on $G$ such that $\pi_n:G\to G_n$ are continuous.

For a topological space $T$, denote by $\widehat{\,T\,}:=\{f: T\to \mathbb S^1\ \mathrm{continuous}\}.$ There is a 
natural compact-open topology on $\widehat{\,T\,}$, generated by open subsets of the form 
$\dis{W(K,U):=\{f\in \widehat{\,T\,}: f(K)\subset U\}}$, where $K\subset T$ are compact, $U\subset\mathbb S^1$ are 
open. We call $\widehat{\,T\,}$  the (topological) dual of $T$.
 
For inductive and projective topologies, we have the following general results concerning their duals.

\begin{prop} 
Let $\{P_n\}_n$ be a projective system of Hausdorff topological groups with structural maps  $\pi_{n,m}: P_n\to P_m$ 
and $\dis{\pi_n:\lim_{\law n}P_n\to P_n}$. Assume that all $\pi_n$ and $\pi_{n,m}$ are surjective and open, and that
for any $n,\,n'$, there exists an $n''$ such that $n''\leq n$ and $n''\leq n'$. Then, as topological groups, 
$$
\widehat{ \lim_{\longleftarrow n} P_n}\simeq \lim_{\longrightarrow n}\widehat{\ P_n\,}.
$$
\end{prop}
\noindent
{\it Proof.} Denote by $\widehat \pi_{n,m}:\widehat{P_m}\to \widehat {P_n},\ f_m\mapsto f_m\circ \pi_{n,m}$,  the dual 
of $\pi_{n,m}: P_n\to P_m$. Then, for an element 
$\dis{\lim_{\longrightarrow n}f_n\in \lim_{\longrightarrow n}\widehat{\ P_n\,}}$, we have $\widehat \pi_{n,m}(f_m)=f_n$, 
or equivalently, $f_m\circ \pi_{n,m}=f_n$. Hence, for an element 
$\dis{x= \lim_{\longleftarrow n} x_n\in \lim_{\longleftarrow n} P_n}$, we have 
$f_m(x_m)=f_m(\pi_{n,m}(x_n))=f_n(x_n)$ for sufficiently small $m,\,n$. Based on this, we define a natural map 
$$
\varphi: \lim_{\longrightarrow n}\widehat{\ P_n\,}\raw \widehat{ \lim_{\longleftarrow n} P_n}, \qquad
\lim_{\longrightarrow n}f_n\mapsto f
$$ 
where $\dis{f:  \lim_{\longleftarrow n} P_n\to\mathbb S^1, \ x= \lim_{\longleftarrow n} x_n\mapsto f_n(x_n).}$  
From the above discussion, $f$ is well defined. Moreover, we have the following

\begin{lem} (1) $f$ is continuous. In particular, $\varphi$ is well defined;

\noindent
(2) $\varphi$ is a bijection;

\noindent
(3) $\varphi$ is continuous; and

\noindent
(4) $\varphi$ is open.
\end{lem}
\noindent
{\it Proof.} (1) Let $U$ be an open subset of $\mathbb S^1$. If $\dis{x= \lim_{\longleftarrow n} x_n\in f^{-1}(U)}$, 
$f_n(x_n)\in U$ for sufficiently small $n$. In particular, $x_n\in f_n^{-1}(U)$. On the other hand, since $f_n$ is 
continuous, $f_n^{-1}(U)$ is open. So, $\dis{\lim_{\longleftarrow n}f_n^{-1}(U)}$ is an open neighborhood of  
$\dis{x= \lim_{\longleftarrow n} x_n}$. Note that $\dis{f(\lim_{\longleftarrow n}f_n^{-1}(U))=f_n(f_n^{-1}(U))=U}$. Hence 
$\dis{\lim_{\longleftarrow n}f_n^{-1}(U)\subset f^{-1}(U)}$. Consequently, $f$ is continuous, and hence 
$\varphi$ is well defined.

(2) To prove that $\varphi$ is injective, we assume that $\dis{\varphi(\lim_{\longrightarrow n}g_n)=:g}=f$. Thus 
$f_n(x_n)=g_n(x_n)$ for sufficiently small $n$ and for all 
$\dis{x= \lim_{\longleftarrow n} x_n\in \lim_{\longleftarrow n} P_n}$. Note that $\pi_{n,m}$ are surjective. So 
$f_n(x_n)=g_n(x_n)$ for all $x_n\in P_n$. This means that $f_n=g_n$ for sufficiently small $n$. Consequently,
$\dis{\lim_{\longrightarrow n}f_n\equiv \lim_{\longrightarrow n}g_n}$, and hence $\varphi$ is injective.

To show that $\varphi$ is surjective,  let $\dis{f:  \lim_{\longleftarrow n} P_n\to\mathbb S^1}$ be a continuous map. 
Then, for any open subset $U\subset\mathbb S^1$ containing 1, $f^{-1}(U)$ is an open neighborhood of 0 in 
$\dis{ \lim_{\longleftarrow n} P_n}$. Hence,  we may write $f^{-1}(U)$ as 
$\dis{f^{-1}(U)=\lim_{\longleftarrow n} P_n\cap\prod_\alpha K_n}$ where $K_n\subset P_n$ are open subsets and 
$K_n=P_n$ for almost all $n$. By assumptions, for $n_1,\dots, n_r$ such that $K_{n_i}=P_{n_i}$, there exists an $N$ such 
that $N\leq n_i$. Then, $f(\mathrm{Ker}\,\pi_N)=1$. So,  $\dis{f( \mathrm{Ker}\,\pi_n)=1}$ for all 
$n\leq N$. Built on this, we define, for $n\geq N$, the maps $f_n:P_n\to\mathbb S^1, x_n\mapsto f(x)$ if $\pi_n(x)=x_n$. Note 
that $f(x)$ always make sense, since $\dis{\pi_n:  \lim_{\longleftarrow n} P_n\to P_n}$ is surjective. Moreover, $f_n$'s are well 
defined. Indeed, if $\dis{y\in  \lim_{\longleftarrow n} P_n}$ such that $\pi_n(y)=x_n$, then $\pi_n(y)=x_n=\pi_n(x)$ for 
$n\leq N$. Hence $x-y\in \mathrm{Ker}\,\pi_n$. This implies that $f(y)=f(x)$. Clearly, by definition, 
$\dis{\varphi(\lim_{\raw n}f_n)=f}$. So $\varphi$ is surjective.

(3) and (4) are direct consequences of the bijectivity of $\varphi$. Indeed, to prove that $\varphi$ is continuous,
it suffices to show that for open subsets of $\dis{\widehat { \lim_{\longleftarrow n} P_n}}$ in the form 
$\frak U=W(K,V)$,  $\varphi^{-1}(\frak U)$ is open in  $\dis{\lim_{\longrightarrow n}\widehat{\ P_n\,}}$, where $K$ is a 
compact subset of $ \dis{\lim_{\longleftarrow n} P_n}$ and $V$ is an open subset of $\mathbb S^1$. Since $\pi_{n}$ are continuous, $K_n=\pi_n(K)$ are 
compact. In this way, we get a inductive system of  open subsets $\{W(K_n,V)\}_n.$ Set 
$\dis{U=\lim_{\raw n} W(K_n,V)}$. Note that, by the bijectivity of $\varphi$, $f(U)= \frak U$. This shows that $f$ is 
continuous.

To prove $\varphi$ is open, let $U$ be an open subset of $\dis{\lim_{\raw n}\widehat{\, P_n\,}}$ such that
 $\big(W(K_n,V)\big)=\iota_n^{-1}(U)$ for a compact subset $K_n$ of $P_n$ for any $n$.  $K:=\lim_{\law_n}K_n$ is compact in 
$\dis{ \lim_{\longleftarrow n} P_n}$. Consequently, $W(K,0)$ is open in $\dis{\widehat { \lim_{\longleftarrow n} P_n}}$. 
Note that, from the bijectivity of $\varphi$, we have $\dis{\varphi(\lim_{\raw n}W(K_n,V))=W(K,V)}$. So $\varphi$ is 
open. This proves the lemma and hence also the proposition.

\vskip 0.25cm
Next we treat inductive systems.  By definition, an inductive system $\{D_n\}_n$ of Hausdorff topological 
groups is called {\it compact oriented}, if for any  compact subset $\dis{K\subset  \lim_{\raw n}{D_n}}$, there 
exists an index $n_0$ such that $K\subset D_{n_0}$.

\begin{prop} 
Let $\{D_n\}_n$ be a compact oriented inductive system of Hausdorff topological groups with  structural  maps  
$\dis{\iota_n:D_n\to \lim_{\raw n}D_n}$  and $\iota_{n, {n'}}:D_n\to D_{n'}$. Assume that  $\iota_{n, {n'}}$  are injective and 
closed, and that, for any $n,\,n'$, there exists an $n''$ such that $n''\geq n$ and $n''\geq n'$. Then, as topological groups, 
$$
\widehat{ \lim_{\raw n} D_n}\simeq \lim_{\law n}\widehat{\ D_n\,}.
$$
\end{prop}
\noindent
{\it Proof.} To start with, we define a map 
$$
\psi: \lim_{\law n}\widehat{\ D_n\,}\raw \widehat{ \lim_{\raw n} D_n},\qquad 
\lim_{\law n}f_n\mapsto f
$$ 
where $\dis{f: \lim_{\raw n} D_n\to\mathbb S^1,\ \lim_{\raw n} x_n\mapsto f_n(x_n).}$

\begin{lem} 
(1) $f$ is well defined and continuous. In particular, $\psi$ is well defined.

\noindent
(2) $\psi$ is a bijection;

\noindent
(3) $\psi$ is continuous; and

\noindent
(4) $\psi$ is open.
\end{lem}
\noindent
{\it Proof.}
(1) Note that for $n<{n'}$, $f_n=f_{n'}\circ\iota_{n,{n'}}$ and $x_{n'}=\iota_{n,{n'}}(x_n)$. Consequently, 
$f_n(x_n)=f_{n'}\circ\iota_{n,{n'}}(x_n)=f_{n'}(\iota_{n,{n'}}(x_n))=f_{n'}(x_{n'})$ for sufficiently large $n\leq n'$. 
So $f$ is well defined. 

To prove that $f$ is continuous, let $U\subset\mathbb S^1$ be an open subset and take 
$\dis{\lim_{\raw n}x_n\in f^{-1}(U)}$. By definition, $x_n\in f_n^{-1}(U)=:U_n$ and $U_n$ are open. 
Hence $\dis{\frak U:=\lim_{\raw n}U_n}$ is open and $\dis{x=\lim_{\raw n}x_n\in\frak U}$. 
Moreover, $f(\frak U)\subset U$. So $f$ s continuous.

(2) Assume that $\dis{\psi(\lim_{\law n}f_n)=f\equiv 0.}$ Then, for any $\dis{x=\lim_{\raw n}x_n}$, 
$f(x)=0$. This means that for all $n$, and $x_n\in D_n$, $f_n(x_n)=f(x)=0$ where $x$ is 
determined by the condition that for all ${n'}\geq n$, $x_{n'}=\iota_{n,{n'}}x_n$. (Since 
$\iota_n$ are injective, this is possible.) Thus $f_n\equiv 0$ and hence $\dis{\lim_{\law n}f_n=0}$.

For any $\dis{f\in \lim_{\law n}\widehat{\ D_n\,}}$, let $f_n=f\circ\iota_n:D_n\to\mathbb S^1$. Clearly, $f_n$ is 
continuous. So $f_n\in\widehat {D_n}$. Moreover, for all ${n'}\geq n$, 
$f_n=f\circ\iota_n=f\circ\iota_{n'}\circ\iota_{n,{n'}}=f_{n'}\circ\iota_{n,{n'}}$. That is, $\{f_n\}_n$ forms a 
projective limit. Obviously, $\dis{\psi(\lim_{\law n}f_n)=f}$.

(3) This is rather involved: not only the just proved bijectivity of $\psi$, but all 
assumptions for our injective system are used here. Let 
$\dis{\frak U=W(K,V)}$ be an open subset of $\dis{\widehat { \lim_{\raw n}{D_n}}}$, where 
$\dis{K\subset \lim_{\raw n}{D_n}}$ is compact and $V$ is an open subset of $\mathbb S^1$. By 
assumptions, for any $n'$, $\iota_{n'}^{-1}(D_n)$ is closed. Hence $\dis{D_n\subset \lim_{\raw n}D_n}$ are closed. 
So $K_n:=K\cap D_n$ are compact. If 
$\dis{\lim_{\law n}f_n}$ is an element in $\dis{\psi^{-1}(\frak U)\subset \widehat{\lim_{\raw n}{D_n}}}$, we have 
$f_n\in W(K_n,V)$, and $\dis{\lim_{\law n}f_n\in \lim_{\law n}W(K_n,V)=\psi^{-1}(W(K,V))}$.  So it suffices to show that  $\dis{\lim_{\law n}W(K_n,V)}$ 
is open. This is a direct consequence of our assumptions. Indeed, 
since our inductive system if compact oriented, there exists a certain $n_0$ such that 
$\dis{K=\lim_{\raw n}K_n\subset D_{n_0}}$. Hence, $K_n=K_{n_0}=K$ for all $n\geq n_0$. 
$\dis{\lim_{\law n}W(K_n,V)=\pi_{n_0}^{-1}(W(K_{n_0},V))}$. Hence, $\dis{\lim_{\law n}W(K_n,V)}$ is open.

(4) This is a direct consequence of the bijectivity of $\psi$. Indeed, let 
$\dis{\frak U\subset \lim_{\law n}\widehat{D_n}}$ be an open subset. By definition, without loss of generality, 
we may assume that there exists an $n$ and an open subset $W(K_n,V)$ of $\widehat{D_n}$ such that 
$\frak U=\pi_n^{-1}(W(K_n,V)).$ Since $K_n$ is a compact subset of $D_n$ and $D_n$ is closed in 
$\dis{\lim_{\raw n}D_n}$, $K_n$ is compact in $\dis{\lim_{\raw n}D_n}$. On the other hand, 
$\psi(\frak U)=\psi(\pi_n^{-1}(W(K_n,V)))=W(K_n,V)$. So $\psi(\frak U)$ is open in $\dis{\widehat{\lim_{\raw n}D_n}}$.
This proves the lemma and hence also the proposition.

\subsubsection{Adelic spaces and their ind-pro topologies}

Let $X$ be an arithmetic surface. For a complete flag $(X, C,x)$ on $X$ (with $C$ an irreducible curve on 
$X$ and $x$ a close point on $C$), let $k(X)_{C,x}$ its associated local ring. By Theorem 1, $k(X)_{C,x}$ is a direct 
sum of two dimensional local fields. Denote by $(\pi_C,t_{x,C})$ a local parameter defined by the flag $C$ of $X$, 
and fix a Madunts-Zhukov lifting ([MZ])
$$
h_C=(h_{\pi_C,t_{X,c}}):\A_{C,01}\simeq{\prod}_{x:x\in C}'\widehat \O_{C,x}\big/\pi_C {\prod}_{x:x\in C}'\widehat \O_{C,x}
\xrightarrow{\mathrm{lifting}}{\prod}_{x:x\in C}'\widehat \O_{C,x}
$$ 
Then, following Parshin, see e.g., Example 2, 
$$
\begin{aligned}
&\A_X^\fin=\A_{X,012}={\prod_{x\in C}}'k(X)_{C,x}:={\prod_C}'\Big({\prod}'_{x:x\in C}k(X)_{C,x}\Big)\\
&:=\Biggl\{\ \big(\sum_{i_C=-\infty}^\infty h_C(a_{i_C})\pi_C^{i_C}\big)_C\in
\prod_C\prod_{x: x\in C}k(X)_{C,x}\,\Big|\,
\begin{aligned}&a_{i_C}\in \A_{C,01},\\ 
&a_{i_C}=0\ (i_C\ll 0);\\ D
&\min\{i_C:a_{i_C}\not=0\}\geq 0\ (\forall' C)
\end{aligned}\ \Biggr\}.
\end{aligned}
$$
This gives the finite adelic space for $X$. Moreover, from Parshin-Osipov, see e.g., Definition 5, we have the infinite 
adelic space $\A_X^\infty:=\A_{X_F}\widehat\otimes_{\mathbb Q}\mathbb R$, and hence the total arithmetic adelic 
space $\A_X^\ar$ for the arithmetic surface $X$:
$$
\A_X^\ar:=\A_X^\fin\oplus\A_X^\infty.
$$
Moreover, there are natural ind-pro structures on $\A_X^\fin$, $\A_X^\infty$ and hence on $\A_X^\ar$, since 
$$
\begin{aligned}
\A_X^\fin=&\lim_{\raw D}\lim_{\law E: E\leq D}\A_{X,12}(D)\big/\A_{X,12}(E),\\
\A_X^\infty=&\lim_{\raw D}\lim_{\law E: E\leq D}\A_{X_F,1}(D)\big/\A_{X_F,1}(E)
\widehat\otimes_{\mathbb Q}\mathbb R.
\end{aligned}$$
Consequently, induced from the locally compact (Hausdorff) topologies on spaces $\A_{X,12}(E)\big/\A_{X,12}(E)$ and 
$\A_{X_F,1}(D)\big/\A_{X_F,1}(E)\widehat\otimes_{\mathbb Q}\mathbb R$, we get canonical ind-pro topologies on 
$\A_X^\fin$, $\A_X^\infty$, and hence on $\A_X^\ar$, which can be easily seen to be Hausdorff. For example, by 
[MZ], a fundamental system of open neighborhood of 0 in $\A_X^\fin $ is given by 
$$
\Biggl\{\ \big(\sum_{i_C=-\infty}^\infty h_C(a_{i_C})\pi_C^{-i_C}\big)_C\in\A_X^\fin :
\begin{aligned}
&a_{i_C}\in U_{i_C}\subset \A_{C,01}\ \mathrm{open\ subgroup}\\
&U_{i_C}=\A_{C,01}\ \forall\ i_C\gg 0\\
& \min\{i_C:U_{i_C}=\A_{C,01}\}\leq 0\ (\forall' C)
\end{aligned}
\ \Biggr\}.\eqno(8)
$$
We have a similar descriptions for $\A_X^\infty$. However, while, with respect to these canonical topologies,  
$\A_X^\fin$, $\A_X^\infty$, and $\A_X^\ar$ are additive topological groups, they are not topological rings.
That is, the multiplication operations are not continuous for these spaces. Still, in \S 3.2.1, we will prove the following 
very useful

\begin{prop} For a fixed ${\bf a}\in\A_C^\ar$, the scalar product by ${\bf a}$, namely, the map
$\A_X^\ar\buildrel{{\bf a}\times}\over\raw\A_X^\ar,\ {\bf x}\mapsto{\bf a}{\bf x}$, is continuous.
\end{prop}

\noindent
{\it Remark.} This result can be used to establish a similar result for two dimensional\\[0.30em] local fields. Indeed,
since $\dis{{\prod}_{x\in C}'k(X)_{C,x}=\lim_{\raw n}\lim_{\law m: m\leq n}\A_{X,12}(nC)\big/\A_{X,12}(mC)}$, as a subspace of
$\dis{\A_X^\fin=\lim_{\raw D}\lim_{\law E: E\leq D}\A_{X,12}(D)\big/\A_{X,12}(E),}$ there is a natural ind-pro topology 
on $k(X)_{C,x}$, induced from the ind-pro topology on $\A_{X,012}.$ Similarly, for a two dimensional
local field $F$, we have $\dis{F=\lim_{\raw n}\lim_{\law m: m\leq n}\frak m_F^{-n}\big/\frak m_F^{-m}}$, where 
$\frak m_F$ denotes the maximal ideal of $F$. So from the natural locally compact topologies on the quotient spaces 
$\frak m_F^{-n}\big/\frak m_F^{-m}$, we obtain yet another ind-pro topology on $F$, and hence on  $k(X)_{C,x}$,
since $k(X)_{C,x}$ is also a direct sum of  two dimensional local fields. Induced from the same roots of locally 
compact topology on one-dimensional local fields, these two topologies on $k(X)_{C,x}$ are equivalent.
This, with the above proposition, then proves the following

\begin{cor} 
For a fixed $a_{C,x}\in k(X)_{C,x}$, the scalar product by $a_{C,x}$, namely, the map
$k(X)_{C,x}\buildrel {a_{C,x}}\cdot\over\longrightarrow k(X)_{C,x},\ \alpha\mapsto a_{C,x}\alpha$  is continuous. In 
particular, the scalar product of a fixed element on a two dimensional local field is continuous.
\end{cor}

We will use this result in an on-going work to prove that, with respect to the canonical ind-pro topology, two dimensional local 
fields are self-dual as topological groups.

\subsubsection{Adelic spaces are complete}

In this section, we show that  adelic spaces $\A_X^\fin $ and $\A_X^\infty$  are complete. For basics of 
complete topological groups, please refer to [Bo] and [G]. We begin with

\begin{prop} The subspaces $\A_{X,12}(D)\subset\A_X^\fin $ and $\A_{X_F,1}(D)\subset\A_X^\infty$, 
and the level two subspace $\A_{X,01}^\ar, \A_{X,02}^\ar$ and $\A_{X,12}^\ar(D)$ of $\A_X^\ar$ 
are complete and hence closed.
\end{prop}
\noindent
{\it Proof.} As our proof below works for all other types as well, we only treat $\A_{X,12}(D)\subset\A_{X}^\fin$ 
to demonstrate. Since $\A_{X,12}(D)\big/\A_{X,12}(E)$'s are finite dimensional vector spaces over one dimensional 
local field, which is locally compact, so they are complete. Consequently, as a projective limit of complete spaces,
$\A_{X,12}(D)=\dis{\lim_{\law E: E\leq D}\A_{X,12}(D)\big/\A_{X,12}(E)}$ is complete. It is also closed 
since $\A_{X,012}$ is Hausdorff.

\begin{prop} 
For an arithmetic surface $X$, its associated adelic spaces $\A_X^\fin $ and $\A_X^\infty$  are complete.
\end{prop}
\noindent
{\it Proof.} We will give a uniform proof for both finite and infinite cases. For this reason, we use simply $\A$ to denote 
both $\A_X^\fin $ and $\A_X^\infty$, and $A(D)$ for both $\A_{X,12}(D)$ and 
$\A_{X_F}(D)\widehat\otimes_{\mathbb Q}\mathbb R$. Clearly, it suffices to prove the following

\begin{lem} 
For a strictly increasing sequence  $\{A(D_n)\}_n$ in $\A$,  $\dis{\lim_{\raw n} A(D_n)}$ is complete. 
\end{lem}
\noindent
{\it Proof.} Let $\{a_n\}_n$ be a Cauchy sequence of $\dis{\lim_{\raw n} A(D_n)}$. We will show that these exists a divisor $D$ 
such that  $\{a_n\}_n\subset A(D)$. Assume that, on the contrary, for all divisors $D$, $\{a_n\}_n\not\subset A(D)$.
We claim that then there exists a subsequence $\{a_{k_n}\}_n$ of $\{a_n\}$, a (strictly 
increasing) subsequence $\{D_{k_n}\}_n$ of $\{D_n\}_n$, and an open  neighborhood  $U$ of 0 in 
$\dis{\lim_{\raw n}A(D_n)}$ such that  (i) $a_{k_n}\in A(D_{k_n})\backslash A(D_{k_{n-1}})$ for all $n\geq 2$; 
(ii) $a_{k_1},\dots, a_{k_n},\dots\not\in U$ and (iii) $a_{k_{i+1}},\dots, a_{k_n},\dots\not\in U+A(D_i), \ i\geq 1$. 
If so, since $a_{k_m}\not\in U+A(D_m)$ and $a_{k_n}\in A(D_m)$, for any $n>m$, $a_{k_n}-a_{k_m}\not\in U$. 
So, $\{a_{k_n}\}_n$ is not a Cauchy sequence of $\dis{\lim_{\raw n} A(D_n)}$. a contradiction. Therefore, there exists a divisor 
$D$ such that $\{a_n\}_n\subset A(D)$. By Proposition 26, $A(D)$ is complete. So the Cauchy sequence $\{a_n\}_n$ is 
convergent in $A(D)$ and hence in $\dis{\lim_{\raw n} A(D_n)}$ as well.

To prove the claim, we select $\{a_{k_n}\}_n$ and the corresponding $D_{k_n}$'s as follows. We begin with
$a_{k_1}=a_1$. Being an element of $\A$, there always a divisor $D_{k_1}$ such that $a_{k_1}\in A(D_{k_1})$.  
Since for all $D$, $\{a_n\}_n\not\subset A(D)$,  there exists $k_2$ and a divisor $D_{k_2}$ such that 
$D_{k_2}>D_{k_1}$ and $a_{k_2}\in A(D_{k_2})-A(D_{k_1})$. By repeating this process, we obtain a subsequence 
$\{a_{k_n}\}_n$ of $\{a_n\}$, a (strictly increasing) subsequence $\{D_{k_n}\}_n$ of $\{D_n\}_n$ such that  (i) above 
holds. Hence to verify the above claim, it suffices to find an open subset $U$ satisfying the conditons (ii) and (iii) 
above. This is the contents of the following

\begin{sublem}
Let $\{A(D_n)\}_n$ be a strictly increasing sequence and $\{a_n\}_n$ be a sequence of elements of $\A$. Assume that
$a_n\in A(D_n)-A(D_{n-1})$ for all $n\geq 1$. Then there exists an open subset $U$ of $\dis{\lim_{\raw n}A(D_n)}$ 
such that $a_1,\dots, a_n,\dots\not\in U$ and $a_{m+1},\dots, a_n,\dots\not\in U+A(D_m)$ for all $m<n$.
\end{sublem}
\noindent
{\it Proof.} We separate the finite and infinite adeles.
\vskip 0.20cm
\noindent
{\bf Finite Adeles\ } Since $A(D_1)/A(D_0)$ is Hausdorff, there exists an open, and hence closed, subgroup 
$U_1\subset A(D_1)$ such that $a_1\not\in U_1$ and $U_1\supset A(D_0)$. Since $A(D_1)$ is complete
and $U_1$ is closed in $A(D_1)$, $U_1$ is complete as well. Now, viewing in $A(D_2)$, since $A(D_2)$ is Hausdorff, 
$U_1$ is a complete subspace, so $U_1$ is  closed in $A(D_2)$. Hence $A(D_2)/U_1$ is Hausdorff too. Therefore,
there exists an open and hence closed subgroup $V_{2,0}$ of $A(D_2)$ such that $a_1, a_2\not\in V_{2,0}$ and 
$V_{2,0}\supset U_1$. In addition, $A(D_2)\big/A(D_1)$ is Hausdorff, there exists an open subgroup $V_{2,1}$ such 
that $a_2\not\in V_{2,1}$ and $V_{2,1}\supset A(D_2)$. Consequently, if we set $U_2=V_{2,0}\cap V_{2,1}$, $U_2$ is 
an open hence closed subgroup of $A(D_2)$ such that $a_1,a_2\not\in U_2$, $a_2\not\in U_2+A(D_1)$ and 
$U_2\supset U_1$. So, inductively, we may assume that there exists an increasing sequence of open subgroups 
$U_1,\dots, U_{n-1}$  satisfying the properties required. In particular,   the following quotient spaces 
$A(D_n)\big/U_{n-1}+A(D_0)\big(=A(D_n)\big/U_{n-1}\big),\dots,A(D_n)\big/U_{n-1}+A(D_m),\dots,
A(D_n)\big/U_{n-1}+A(D_{n-1})\big(=A(D_n)\big/A(D_{n-1})\big)$ 
are Hausdorff. Hence, there are open subgroups $V_{n,m},\ 0\leq m\leq n-1$ of $A(D_n)$ such that 
$a_{m+1},\dots, a_n\not\in V_{n,m}$ and $V_{n,m}\supset U_{n-1}+A(D_m)$. Define 
$U_n:=\cap_{m=1}^{n-1}V_{n,m}$. Then $U_n$ is an open subgroup of $A(D_n)$ satisfying 
$a_1,\dots, a_n\not\in U_n$, $a_{m+1},\dots, a_n\not\in U_n+A(D_m), \ 1\leq m\leq n-1$ and $U_n\supset U_{n-1}$.
Accordingly, if we let $\dis{U=\lim_{\raw n}U_n}$, by definition, $U$ is  an open subgroup of 
$\dis{\lim_{\raw n}A(D_n)}$, and  from our construction, $a_1, \dots, a_n,\dots\not\in U$ and 
$a_{m+1},\dots, a_n,\dots\not\in U+A(D_m),\ m\geq 1$.
\vskip 0.20cm
\noindent
{\bf Infinite Adeles\ } Since $A(D_1)$ is Hausdorff, there exists an open subset $U_1$ of $A(D_1)$ such that 
$a_1\not\in A(D_1)$. Moreover, since $A(D_2)\simeq A(D_2)\big/A(D_1)\oplus A(D_1)$ and $A(D_2)\big/A(D_1)$ is 
Hausdorff, there exists an open subset $U_2$ of $A(D_2)$ such that $a_1, a_2\not\in U_2$ and 
$U_2\cap A(D_1)=U_1$. In particular, $a_2\not\in U_2+A(D_1)$. Similarly, as above, with an inductive process, based 
on the fact that $A(D_n)\simeq A(D_n)\big/A(D_{n-1})\oplus A(D_{n-1})$ and $A(D_n)\big/A(D_{n-1})$ is Hausdorff, 
there exists an open subset $U_n$ of $A(D_n)$ such that $a_1,\dots,a_n\not\in U_n$ and $U_n\cap A(D_n)=U_{n-1}$. 
Consequently, $a_{m+1},\dots, a_n\not\in U_n+A(D_m),\ 1\leq m\leq n-1$. In this way, we obtain an infinite increasing 
sequence of open subsets $U_n$. Let $\dis{U=\lim_{\raw n}U_n}$. Then by definition $U$ is an open subset of 
$\dis{\lim_{\raw n}A(D_n)}$ satisfying $a_1,\dots, a_n,\dots\not\in U$ and 
$a_{m+1},\dots, a_n,\dots\not\in U+A(D_m),\ m\geq 1$. This then proves the sublemma, the lemma and hence also 
the proposition.

\subsubsection{Adelic spaces are compact oriented}

In this section, we show that adelic spaces $\A_X^\fin $ and $\A_X^\infty$  are compact oriented.

\begin{prop} 
For an arithmetic surface $X$, its associated adelic spaces $\A_X^\fin $ and $\A_X^\infty$  are compact oriented. 
That is to say, for any   compact subgroup, resp. a compact subset, $K$ in $\A_X^\fin $, resp. $\A_X^\infty$, there 
exists a divisor $D$ on $X$, resp., on $X_F$, such that $K\subset \A_{X,12}(D)$, resp., $K\subset \A_{X_F}(D)$.
\end{prop}
\noindent
{\it Proof.} We treat both finite and infinite places simultaneously. Assume that for all divisors 
$D$, $K\not\subset A(D)$. Fix a suitable $D_0$. By our assumption, $K\not\subset A(D_0)$. Since $\A$ is Hausdorff 
and $A(D_0)$ is closed, there exists a certain divisor $D_1$ and an element $a_1\in K$, such that  $D_1>D_0$, 
$a_1\in A(D_1)\backslash A(D_0)$. Similarly, since $A(D_1)$ is closed, $\A\big/A(D_0)$  is Hausdorff, we can find an 
open subgroup $U_1'\subset \A\big/A(D_0)$ such that $a_1+A(D_0)\not\subset \A\big/A(D_0)$. Consequently, there 
exists an open subgroup $U_1$ of $\A$ such that $a_1\not\in U_1$ and $U_1\supset A(D_0)$. Now use $(D_1,\A)$ 
instead of $(D_0,\A)$, by repeating the above construction, we can find a divisor $D_2$, an open subgroup 
$U_2\subset\A$ and an element $a_2\in K\cap \big(A(D_2)\backslash A(D_1)\big)$ such that $D_2>D_1$,
$a_2\not\in U_2$, $U_2\supset U_1$. In this way, we obtain a sequence of divisors $D_n$, a sequence of elements 
$a_n\in K\cap \big(A(D_n)\backslash A(D_{n-1})\big)$ and a sequence of open subgroups $U_n$ such that 
$D_\alpha>D_{n-1}$,  $a_n\not\in U_n$, $U_n\supset U_{n-1}$. Let $U=\dis{\lim_{\raw n}U_n}$. Then $U$ is an 
open, and hence closed, subgroup of $\A$. Consequently, $K\cap U$ is  compact. This is a contradiction. Indeed, 
since $a_1,\dots, a_n\not\in U_n$ for all $n$, the open covering $\{U_n\}_n$ of $K\cap U$ admits no finite 
sub-covering. This completes the proof.

\subsubsection{Double dual of adelic spaces}

We here prove the following

\begin{prop} 
As topological groups, we have
$$\begin{aligned}
\Big(\lim_{\law E: E\leq D}\A_{X,12}( D)\big/\A_{X,12}( E)\Big)^\vee &\simeq\  
\lim_{\raw E: E\leq D}\Big(\A_{X,12}( D)\big/\A_{X,12}( E)\Big)^\vee,\\
\Big(\lim_{\law E: E\leq D}\big(\A_{X_F}(D)\big/\A_{X_F}(E)\widehat\otimes_{\mathbb Q}\mathbb R\big)\Big)^\vee&\simeq\  
\lim_{\raw E: E\leq D}\Big(\A_{X_F}(D)\big/\A_{X_F}(E)\widehat\otimes_{\mathbb Q}\mathbb R\Big)^\vee;\end{aligned}
$$ 
and
$$\begin{aligned}
\Big(\lim_{\raw  D}\A_{X,12}( D)\Big)^\vee& \simeq\ \lim_{\law D}\Big(\A_{X,12}( D)\Big)^\vee,\\
\Big(\lim_{\raw  D}\big(\A_{X_F}(D)\widehat\otimes_{\mathbb Q}\mathbb R\big)\Big)^\vee
&\simeq\ \lim_{\law D}\Big(\A_{X_F}(D)\widehat\otimes_{\mathbb Q}\mathbb R\Big)^\vee.\end{aligned}
$$
where, $\dis{\A_{X_F}(D)\widehat\otimes_{\mathbb Q}\mathbb R
:=\lim_{\law E: E\leq D}\Big(\A_{X_F}(D)\big/\A_{X_F}(E)\widehat\otimes_{\mathbb Q}\mathbb R}\Big)$.

In particular, 
$$\begin{aligned}
\widehat{\A_X^\fin }&\simeq\lim_{\law D}\lim_{\substack{\raw E\\ E\leq D}}\Big(\A_{X,12}( D)\big/\A_{X,12}( E)\Big)^\vee,\\ 
\widehat{\A_X^\infty}&\simeq\lim_{\law D}\lim_{\substack{\raw E\\ E\leq D}}
\Big(\A_{X_F}(D)\big/\A_{X_F}(E)\widehat\otimes_{\mathbb Q}\mathbb R\Big)^\vee.\end{aligned}
$$
\end{prop}
\noindent
{\it Proof.} We apply Proposition 20, resp. Proposition 22, to prove the first, resp., the second, pairs of  
homeomorphisms. We need to check the conditions  there.

As above, we treat finite adeles and infinite adeles simultaneously. So we use $\A$ and $A(D)$ as in
\S 3.1.3. Then $\dis{A(D)=\lim_{\law E: E\leq D}A(D)\big/A(E)}$.  Now, for $E<E'$, 
$\pi_{D/E,D/E'}:A(D)\big/A(E)\to A(D)\big/A(E')$ is the natural quotient map. So  $\pi_{D/E,D/E'}$ are surjective and open. 
Similarly,  $\pi_{D,D/E}:A(D)\to A(D)\big/A(E)$ are surjective and open. So, by Proposition 20, we get the first group of two 
homeomorphisms for topological groups.

To treat the second group, recall that $A(D)\big/A(E)$ are complete. So, their projective limits $A(D)$'s are complete.
This implies that $\A(D)$ are closed in $\A$, since $\A$ is Hausdorff. On the other hand, for $D<D'$, 
$\A_{X,12}( D)\subset \A_{X,12}( D')$. So the structural maps $\iota_{D,D'}: \A(D)\to\A( D')$ and $\iota: A(D)\to \A$
are injective and closed. Thus, by Proposition 22, it suffices to show that the inductive system $\{A(D)\}_D$ is compact 
oriented. This is simply the contents of \S3.1.3-4. All this then completes our proof, since the last two homeomorphisms 
are direct consequences of previous four.

\begin{cor} As topological groups,
$$
\widehat{\ \widehat{\A_{X}^\fin}\ }\simeq \A_{X}^\fin,\qquad 
\widehat{\ \widehat{\ \A_{X}^\infty\ }\ }\simeq \A_{X}^\infty, \qquad\mathrm{and\ hence}\qquad \widehat{\ 
\widehat{\A_{X}^\ar}\ }\simeq \A_{X}^\ar.
$$
\end{cor}
\noindent
{\it Proof.} Since $A(D)/A(E)$ are locally compact and hence they are self dual. Thus, to prove this double dual 
properties for our spaces, it suffices to check the conditions listed in Proposition 22 for inductive systems
$\dis{\{\lim_{\raw E:E\leq D}\widehat{A(D)\big/A(E)}\}_E}$ and in Proposition 20 for the projective system 
$\{\widehat {A(D)}\}_D$. With the above lengthy discussions, all this now becomes rather routine. For example, to 
verify that $\widehat {\,\A\,}$ is Hausdorff, we only need to recall that $\mathbb S^1$ is compact.  Still, as careful 
examinations would help understand the essences of our proof above, we suggest ambitious readers to supply 
omitted details.

\subsection{Adelic spaces and their duals}

\subsubsection{Continuity of scalar products}

We here show that Proposition 24, namely, the scalar product maps on adelic spaces are continuous, even adelic spaces are 
not topological rings.
\vskip 0.25cm
\noindent
{\bf Proposition 24.} {\it For a fixed element ${\bf a}$ of $\A_X^\fin$, resp. of $\A_X^\infty$, the induced scalar product map: 
$\phi_{\bf a}^\fin:\A_X^\fin \buildrel {\bf a}\times\over\longrightarrow \A_X^\fin$, resp., 
$\phi_{\bf a}^\infty:\A_X^\infty \buildrel {\bf a}\times\over\longrightarrow \A_X^\infty$, is continuous.}
\vskip 0.25cm
\noindent
{\it Proof.} If ${\bf a}=0$, there is nothing to prove.  Assume, from now on, that 
${\bf a}=(a_C)_C\dis{=(\sum_{i_C=i_{C,0}}^\infty  h_C(a_{i_C})\pi_C^{i_C})_C\in \A_X^\fin }\not=0$. Here, for each 
$C$, we assume that $a_{i_{C,0}}\not=0$. To prove that $\varphi_a$ is continuous, by our description (8) of the 
ind-pro topology on $\A_X^\fin $, it suffices to show that for an open subgroup 
$U=(U_C)_C\dis{=\big(\sum_{j_C=-\infty}^\infty h_C(\A_{C,1}(D_{j_C})\big)\pi_C^{i_C}
+\sum_{j_C=r_C}^\infty h_C(\A_{C,01})\pi_C^{i_C}\big)\cap\A_X^\fin }$, as an open neighborhood of 0, 
its inverse image $\varphi_{\bf a}^{-1}(U)$ contains an open subgroup. For later use, set $I_C:=r_C-i_{C,0}$.

Let $\dis{{\bf b}=(b_C)_C=\big(\sum_{k_C=-\infty}^\infty h_C(b_{k_C})\pi_C^{k_C}\big)_C
\in \varphi_{\bf a}^{-1}(U)\subset\A_X^\fin.}$ Then, for each fixed $C$,  
$\dis{a_Cb_C=\sum_{l_C=-\infty}^\infty\big(\sum_{i_C=i_{C,0}}^\infty h_C(a_{i_C})h_C(b_{l_C-i_C})\big)\pi_C^{l_C}}$. 
Recall that $h_C$ is the lifting map
$\dis{h_C:\A_{C,01}\simeq{\prod_{x:x\in C}}'\widehat \O_{C,x}\big/\pi_C {\prod_{x:x\in C}}'\widehat \O_{C,x}
\xrightarrow{\mathrm{lifting}}{\prod_{x:x\in C}}'\widehat \O_{C,x}.}$ Thus if $b_{k_C}\in \A_{C,01}$, we always have
$h_C(a_{i_C})h_C(b_{l_C-i_C})\in \dis{\sum_{m_C=0}^\infty h_C\big(\A_{C,01}\big)\pi_C^{m_C}}$. 
Moreover, if we write, as we can,  $a_{i_C}\in \A_{C,1}(F_{i_C}),\ b_{k_C}\in \A_{C,1}(E_{k_C})$ for some divisors 
$F_{i_C}$ and $E_{k_C}$, we have $\dis{h_C(a_{i_C})h_C(b_{l_C-i_C})\in 
\sum_{m_C=0}^\infty h_C\big(\A_{C,1}(F_{i_C}+E_{l_C-i_C})\big)\pi_C^{m_C}.}$

Now write $b_C=\big(\sum^{I_C-1}_{k_C=-\infty}+\sum_{k_C= I_C}^\infty\big) h_C(b_{k_C})\pi_C^{k_C}$. We will 
construct the required open subgroup according to the range of the degree index $k_C$.

\noindent
(i) If $\dis{b_C\in \big(\sum_{k_C=I_C}^\infty h_C(\A_{C,01})\pi_C^{k_C}\big)
\cap\big({\prod_{x:x\in C}}'k(X)_{C,x}\big)}$,  we have $a_Cb_C\in U_C$;

\noindent
(ii) To extend the range including also the degree $I_C-1$,  choose a divisor $E_{I_C-1}$ such that
$h_C\big(\A_{C,1}(F_{i_{C,0}}+E_{I_C-1})\big)\subset h_C\big(\A_{C,1}(D_{r_C-1})\big)$. Then, if we choose
$b_C\in\dis{\big(h_C\big(\A_{C,1}(E_{I_C-1})\pi_C^{I_C-1}
+\sum_{k_C=I_C}^\infty h_C(\A_{C,01})\pi_C^{k_C}\big)\cap\big({\prod_{x:x\in C}}'k(X)_{C,x}\big)}$,
we also have $a_Cb_C\in U_C$;

\noindent
(iii) Similarly, to extend the range including the degree $I_C-2$, choose a divisor $E_{I_C-2}$ such that 
$h_C\big(\A_{C,1}(F_{i_{C,0}}+E_{I_C-2})\big)\subset h_C\big(\A_{C,1}(D_{r_C-2})\big)\cap h_C\big(\A_{C,1}
(D_{r_C-1})\big)$ and $h_C\big(\A_{C,1}(F_{i_{C,0}+1}+E_{I_C-2})\big)\subset h_C\big(\A_{C,1}(D_{r_C-1})\big)$. 
Then,  if we choose $b_C\in\dis{\big(\sum_{k_C=I_C-2}^{I_C-1}h_C\big(\A_{C,1}(E_{k_C})\pi_C^{k_C}
+\sum_{k_C=I_C}^\infty h_C(\A_{C,01})\pi_C^{k_C}\big)\cap\big({\prod_{x:x\in C}}'k(X)_{C,x}\big)}$, then we have 
$a_Cb_C\in U_C$.

\noindent
Continuing this process repeatedly, we then obtain divisors $E_{k_C}$'s such that, for 
$b_C\in V_C\dis{:=\big(\sum_{k_C=-\infty}^{I_C-1}h_C\big(\A_{C,1}(E_{k_C})\pi_C^{k_C}
+\sum_{k_C=I_C}^\infty h_C(\A_{C,01})\pi_C^{k_C}\big)\cap\big({\prod_{x:x\in C}}'k(X)_{C,x}\big)}$,  we have 
$a_Cb_C\in U_C$. 

Since, for all but finitely many $C$, $r_C\leq 0$ and $i_{C,0}\geq 0$, or better, $I_C=r_C-i_{C,0}\leq 0$. Therefore, 
from above discussions, we conclude that $\prod_CV_C\cap\A_X^\fin $ is an open subgroup of $\A_X^\fin $ and
$a\big(\prod_CV_C\cap\A_X^\fin \big)\subset U$. In particular, $\phi_{\bf a}$ is continuous.

A similar proof works for $\phi_{\bf a}^\infty$. We leave details to the reader. 

\subsubsection{Residue maps are continuous}

Fix  a non-zero rational differential $\omega$ on $X$. Then for an element ${\bf a}$ of $\A_X^\fin$, resp., $\A_X^\infty$,  
induced from the natural residue pairing $\langle\cdot,\cdot\rangle_\omega$, we get a natural map
$\varphi_{\bf a}^\fin:=\langle{\bf a},\cdot\rangle_\omega:\A_X^\fin \raw\mathbb R/\Z$, resp., 
$\varphi_{\bf a}^\infty:=\langle{\bf a},\cdot\rangle_\omega:\A_X^\infty \raw\mathbb R/\Z$.

\begin{lem}
Let ${\bf a}$ be a fix element in $\A_X^\fin$, resp., $\A_X^\infty$. Then the induced map
$\varphi_{\bf a}^\fin:=\langle{\bf a},\cdot\rangle_\omega:\A_X^\fin \raw\mathbb R/\Z$, resp., 
$\varphi_{\bf a}^\infty:=\langle{\bf a},\cdot\rangle_\omega:\A_X^\infty \raw\mathbb R/\Z$,  is 
continuous. In particular, the residue map on arithmetic adeles $\A_X^\ar$ is continuous.
\end{lem}
\noindent
{\it Proof.} We prove only for $\varphi_{\bf a}^\fin$, as a similar proof works $\varphi_{\bf a}^\infty$.
Write $\A_X^\fin =\prod_{\substack{F:\,\mathrm{2-dim}\\ \ \mathrm{local\ field}}}'F$. And, for each local field  $F$, fix 
an element $t_F$ of $F$ such that for equal characteristic field $F$,  $t_F$ is a uniformizer of $F$, while for 
mixed characteristics field $F$, $t_F$ is a lift of a uniformlzer of its residue field. Since, by Proposition 27, the scalar 
product is continuous,  to prove the continuity of $\langle{\bf a},\cdot\rangle_\omega$,  it suffices to show that the  
residue map $\Res:\A_X^\fin \to\mathbb R/\Z,\ {\bf x}=(x_F)\mapsto\sum_F\res_F(x_F\,dt_F)$ is continuous.
(Note that, by the definition of $\A_X^\fin$, see, e.g., \S 3.1.2, the above summation is a finite sum.)  Since the open subgroup 
$\Big(\sum_{i_C=\infty}^{-1}h_C(\A_{C,1}(0)\big)\pi_C^{i_C}+\sum_{i=0}^\infty h_C(\A_{C,01})\pi_C^{i_C}\Big)
\cap\A_X^\fin $ is contained, the kernel of the residue map is an open subgroup. This proves the lemma.

\subsubsection{Adelic spaces are self-dual}

We will treat both $\A_X^\fin $ and $\A_X^\infty$ simultaneously. So as before, we use $\A$ to represent them.

Recall that, for a fixed ${\bf a}\in \A$, the map $\langle {\bf a},\cdot\rangle_\omega:\A\to\mathbb S^1$ is continuous.
Accordingly, we define a map 
$\varphi:\A\to \widehat{\A},\ {\bf a}\mapsto  \varphi_{\bf a}:=\langle {\bf a},\cdot\rangle_\omega.$ 

\begin{prop} 
For the map $\varphi:\A\to \widehat{\A},\ {\bf a}\mapsto  \varphi_{\bf a}:=\langle {\bf a},\cdot\rangle_\omega,$ we have

\noindent
(1) $\varphi$ is continuous;

\noindent
(2) $\varphi$ is injective;

\noindent
(3) The image of $\varphi$ is dense;

\noindent
(4) $\varphi$ is open.
\end{prop}
\noindent
{\it Proof.} (1) For an open subset $W(K,0)$ of $\widehat{\,\A\,}$, where $K$ is a compact subgroup, resp. a compact subset, 
of $\A$, let $U:=\varphi^{-1}\big(W(K,0)\big)$. By Proposition 31,  
$\dis{\widehat {\,\A\,}=\lim_{\law D}\lim_{\raw E: E\leq D}\widehat{A(D)/A(E)}}$. So we may write 
$\chi_0:=\langle {\bf 1},\cdot\rangle_\omega$ as $\dis{\lim_{\law D}\lim_{\raw E: E\leq D}\chi_{D/E}}$ with 
$\chi_{D/E}\in \widehat{A(D)/A(E)}$. Accordingly, write $A_{D/E}:=A(D)/A(E)$, $K_{D/E}:=K\cap A(D)\big/ K\cap A(E)$ 
and let 
$U_{D/E}:=\big\{a_{D/E}\in A_{D/E}\,:\,\chi_{D/E}\big(a_{D/E}K_{D/E}\big)
=\{0\}\ \mathrm{resp.\ an\ open\ subset}\ V\big\}$. Since, for a fixed divisor $D$, $A(D)$ is closed in $\A$,  
$K\cap A(D)$ is a  subgroup, resp. a subset, of $A(D)$. So, for $E\leq D$, $K_{D/E}$ is compact in $A_{D/E}$. 
Consequently, from the non-degeneracy of $\chi_{D/E}$ on locally compact spaces, $U_{D/E}$ is an open subgroup, 
resp., an open subset, of $\A$, and $\dis{U=\lim_{\law D}\lim_{\raw E: E\leq D}U_{D/E}}$. We claim that $U$ is open. 
Indeed, by Proposition 30, $\A$ is compact oriented. So, for  compact $K$, there exists a divisor $D_1$ such that 
$K\subset A(D_1)$. On the other hand, since $\chi_0$ is continuous, there exists a divisor $D_2$ such that 
$A(D_1+D_2)\subset \mathrm{Ker}(\chi_0)$. Hence $U\supset A(D_2)$. Thus, for a fixed $D$, with respect to 
sufficiently small $E\leq D$, we have $U_{D/E}=A_{D/E}$. This verifies that $U$ is open, and hence proves (1), since 
the topology of $\widehat{\,\A\,}$ is generated by the open subsets of the form $W(K,0)$.
\vskip 0.20cm
(2) is a direct consequence of the non-degeneracy of the residue pairing. So we have (2). 

To prove (3), we use the fact that $\A\buildrel\psi\over\simeq\widehat{\ \widehat{\,\A\,}\ }$, where, for ${\bf a}\in\A$, 
$\psi_{\bf a}$ is given by $\psi_{\bf a}:\widehat{\,\A\,}\to\mathbb S^1, \ \chi\mapsto\chi({\bf a})$. Thus to show that the 
image of $\varphi$ is dense, it suffices to show that the annihilator subgroup 
$\mathrm {Ann}\big(\mathrm{Im}(\varphi)\big)$ of $\mathrm{Im}(\varphi)$ is zero. Let then 
${\bf x}\in \mathrm {Ann}\big(\mathrm{Im}(\varphi)\big)$ be an annihilator of $\mathrm{Im}(\varphi)$. Then, by 
definition, $\{0\}=\psi_{\bf x}\big(\{\varphi_{\bf a}:{\bf a}\in\A\}\big)=\{\varphi_{\bf a}({\bf x}):{\bf a}\in\A\}$. That is to say, 
$\langle {\bf a},{\bf x}\rangle_\omega=0$ for all ${\bf a}\in\A\}$. But the residue pairing is non-degenerate. So, 
${\bf x}=0$. 
\vskip 0.20cm
(4) This is the dual of (2). Indeed, let $U\subset \A$ be an open subgroup, resp. an open subset, of $\A$.
Then $U\cap A(D)$ is open in $A(D)$. Since $A(D)$ is closed, $U_{D/E}:=U\cap A(D)\big/ U\cap A(E)$ is open in 
$A_{D/E}$. This, together with the fact that $\chi_{D/E}$ is non-degenerate on its locally compact base space, implies 
that $K_{D/E}:=\big\{a_{D/E}\in A_{D/E}:\chi_{D/E}\big(a_{D/E}\cdot U_{D/E}\big)
=\{0\}\ \mathrm{resp.\ an\ open\ subset}\ V \big\}$ is a compact subset, resp. a compact subset. 
Let $\dis{K:=\lim_{\raw D}\lim_{\law E: E\leq D}K_{D/E}}$. Since $U$ is open, there exists a divisor $E$ such that 
$A(E)\subset U$. This implies that there exists a divisor $D$ such that $K\subset A(D)$. Otherwise,  assume that, for 
any $D$, $K\not\subset A(D)$. Then, there exists an element $k\in K$ such that  $k\not\in A (\omega)-E)$. Hence we 
have $\chi(k\cdot A(E))\not=\{0\}$, a contradiction. This then completes the proof of (4), and hence the proposition.
\vskip 0.20cm
We end this long discussions on ind-pro topology over adelic space $\A_X^\ar$ with the following main theorem.

\begin{thm} Let $X$ be an arithmetic surface. Then, as topological groups, 
$$
\widehat{\ \A_X^\fin \ }\simeq \A_X^\fin 
\quad\mathrm{and}\quad 
\widehat{\ \A_X^\infty\ }\simeq \A_X^\infty.
\qquad{In\ particular,}\quad
\widehat{\ \A_X^\ar\ }\simeq \A_X^\ar.
$$
\end{thm}
\noindent
{\it Proof.} With all the preparations above, this now becomes rather direct. Indeed, by Proposition 34, 
we have an injective continuous open morphism $\varphi:\A\to\widehat{\,\A\,}.$ So it suffices to
show that $\varphi$ is surjective. But this is a direct consequence of the fact that $\varphi$ is dense, 
since both $\A$ and $\widehat{\,\A\,}$ are complete and Hausdorff. This proves the theorem and hence also Theorem II.

\subsubsection{Proof of cohomological duality}

Now we are ready to prove Theorem 19, or the same Theorem III, for the duality of cohomology groups. Recall that for a 
non-zero rational differential $\omega$, by  Proposition 12, we have a non-degenerate residue pairing 
$$
\langle\cdot,\cdot\rangle_\omega:  \A_X^\ar\times  \A_X^\ar\raw\mathbb S^1.
$$
Moreover, by Theorem II just proved, we obtain a natural homeomorphism of topological groups
$$
\A_X^\ar\ \simeq\ \widehat{\ \A_X^\ar\ },\quad {\bf a}\mapsto\langle {\bf a},\cdot \rangle_\omega.
$$
This,  with a well-known argument which we omit, implies the following

\begin{lem} 
With respect to the non-degenerate pairing $\langle\cdot,\cdot\rangle$ on $\A_X^\ar$, we have

\noindent
(i) If $\,W_1$ and $W_2$ are closed subgroups of $\A_X^\ar$, 
$$
(W_1+W_2)^\perp=W_1^\perp\cap W_2^\perp\qquad\mathrm{ and}\qquad
(W_1\cap W_2)^\perp=W_1^\perp+ W_2^\perp;
$$

\noindent
(ii) If  $\,W$ is a closed subgroup of $\A_X^\ar$, then, algebraically and topologically,
$$
(W^\perp)^\perp=W\qquad\mathrm{ and}\qquad W\simeq \,\widehat {\ \A_X^\ar\big/ W^\perp\ }.
$$
\end{lem}

With this, we can complete our proof,  using Proposition 15 for perpendicular subspaces of our level two subspaces 
$\A_{X,01}^\ar,\, \A_{X,02}^\ar, \,\A_{X,12}^\ar(D)$ of $\A_X^\ar$, as follows:

\noindent
(1) {\bf Topological duality between $H_\ar^0$ and $H_\ar^2$}
$$
\begin{aligned}
\widehat {H^2_{\ar}(X,(\omega)-D)}\simeq& \Big(\A_{X,01}^\ar+\A_{X,02}^\ar+\A_{X,12}^\ar((\omega)-D)\Big)^\perp\\
\simeq &\Big(\A_{X,01}^\ar\Big)^\perp\cap\Big(\A_{X,02}^\ar\Big)^\perp\cap\Big(\A_{X,12}^\ar((\omega)-D)\Big)^\perp\\
=&\A_{X,01}^\ar\cap\A_{X,02}^\ar\cap \A_{X,12}^\ar(D)\simeq H^{0}_{\ar}(D);
\end{aligned}
$$
 
\noindent
(2) {\bf Topological duality among $H_\ar^1$}
$$
\begin{aligned}
\widehat{H^1_{\ar}(X,(\omega)-D)}
=&\Big(\frac{\A_{X, 02}^{\ar}\cap \big(\A_{X, 01}^{\ar}+\A_{X, 12}^{\ar}((\omega)-D)\big)}
{\A_{X, 01}^{\ar}\cap \A_{X, 02}^{\ar}+\A_{X, 02}^{\ar}\cap \A_{X, 12}^{\ar}((\omega)-D)}\Big)^{\widehat{~~~}}\\
\simeq&\frac{\big(\A_{X, 01}^{\ar}\cap \A_{X, 02}^{\ar}\big)^\perp\cap
\big(\A_{X, 02}^{\ar}\cap \A_{X, 12}^{\ar}((\omega)-D)\big)^\perp}
{\big(\A_{X, 02}^{\ar}\big)^\perp+ \big(\A_{X, 01}^{\ar}+\A_{X, 12}^{\ar}((\omega)-D)\big)^\perp}\\
=&\frac{\big(\A_{X, 01}^{\ar}+ \A_{X, 02}^{\ar}\big)\cap\big(\A_{X, 02}^{\ar}+ \A_{X, 12}^{\ar}(D)\big)}
{\A_{X, 02}^{\ar}+\A_{X, 01}^{\ar}\cap\A_{X, 12}^{\ar}(D)}\\
\simeq&\frac{\big(\A_{X, 01}^{\ar}+\A_{X, 02}^{\ar}\big)\cap \A_{X, 12}^{\ar}(D)}
{\A_{X, 01}^{\ar}\cap \A_{X, 12}^{\ar}(D)+\A_{X, 02}^{\ar}\cap \A_{X, 12}^{\ar}(D)}\simeq H^1_{\ar}(X,D).
\end{aligned}$$ 
This then completes the proof of Theorem III.
\vskip 0.20cm
{\footnotesize{
\noindent
{\bf Acknowledgements.} 
Special thanks due to Osipov for explaining his joint works with Parshin. LW would like to 
thank  Hida and Mathematics Department, UCLA, for providing him an excellent working environment during the 
preparation for the first version of this work. 

This work is partially supported by JSPS.
}}

\eject
\noindent
{\bf REFERENCES}
\vskip 0.20cm
\noindent
[B] A.A. Beilinson, Residues and adeles, Funct. Anal. Pril., 14 (1980), no. 1, 44-45; English transl. in 
Func. Anal. Appl., 14, no. 1 (1980), 34-35.
\vskip 0.10cm
\noindent
[Bo] N. Bourbaki, {\it Elements of Mathematics: General Topology}, Springer, 1987
\vskip 0.10cm
\noindent
[G] S.A. Gaal, {\it Point set topology}, Pure and Applied Mathematics, Vol. XVI, Academic Press,1964
\vskip 0.10cm
\noindent
[H] A. Huber,  On the Parshin-Beilinson adeles for schemes, Abh. Math. Sem. Univ. Hamburg, 61 (1991), 249-273.
\vskip 0.10cm
\noindent
[Iw] K. Iwasawa, {\it Algebraic Theory of Algebraic Functions,} Lecture notes at Princeton university, 1975
\vskip 0.10cm
\noindent
[L] S. Lang, {\it Introduction to Arakelov theory},  Springer Verlag, 1988
\vskip 0.10cm
\noindent
[MZ] A. I. Madunts and I. B. Zhukov, Multidimensional complete fields: topology and other
basic constructions, Trudy S.-Peterb. Mat. Obshch. (1995); English translation in Amer.
Math. Soc. Transl. (Ser. 2) 166 (1995), 1-34.
\vskip 0.10cm
\noindent
[M1] M. Morrow, An explicit approach to residues on and dualizing sheaves of arithmetic surfaces. New York J. Math. 16 (2010), 575-627.
\vskip 0.10cm
\noindent
[M2] M. Morrow, Grothendieck's trace map for arithmetic surfaces via residues and higher adeles, Journal of Algebra and 
Number Theory,  6 (2012) 1503-1536
\vskip 0.10cm
\noindent
[OP] D. V. Osipov, A. N. Parshin, Harmonic analysis on local fields and adelic spaces. II.  Izv. Ross. Akad. 
Nauk Ser. Mat. 75 (2011), no. 4, 91--164
\vskip 0.10cm
\noindent
[P1] A.N. Parshin,  On the arithmetic of two-dimensional schemes. I. Distributions and residues.  Izv. Akad. 
Nauk SSSR Ser. Mat. 40 (1976), no. 4, 736-773
\vskip 0.10cm
\noindent
[P2] A.N. Parshin, Chern classes, adeles and $L$-functions. J. Reine Angew. Math. 341 (1983), 174-192.
\vskip 0.10cm
\noindent
[S] J.P. Serre, Algebraic groups and class fields, GTM 117, 1988, Springer-Verlag
\vskip 0.10cm
\noindent
[SW] K. Sugahara and L. Weng, Topological duality of two dimensional local fields, in preparations
\vskip 0.10cm
\noindent
[T] J. Tate,  Fourier analysis in number fields, and Hecke's zeta-functions. 1967 Algebraic Number Theory 
pp. 305-347 Thompson, Washington, D.C.
\vskip 0.10cm
\noindent
[W] L. Weng, Geometry of numbers,  arXiv:1102.1302
\vskip 0.10cm
\noindent
[Y] A. Yekutieli, {\it An Explicit construction of the Grothendieck residue complex}, with an appendix by 
P. Sastry,  Asterisque No. 208 (1992)
\vskip 0.50cm
K. Sugahara and L. Weng

Graduate School of Mathematics,
 
Kyushu University,
 
Fukuoka, 819-0395,
 
Japan
 
E-mails: k-sugahara@math.kyushu-u.ac.jp,
 
\hskip 2.30cm  weng@math.kyushu-u.ac.jp
\end{document}